\documentclass[a4paper,fleqn]{cas-sc}

\usepackage[numbers]{natbib}
\usepackage{xcolor}
\usepackage{pdflscape}
\usepackage{amsmath,amssymb,amsfonts,amsthm}
\usepackage{graphicx}
\usepackage{float}
\usepackage{xcolor}
\usepackage{hyperref}
 \hypersetup{colorlinks,citecolor=blue,linkcolor=blue,breaklinks=true}
\usepackage{booktabs}
\usepackage{colortbl}
\usepackage{enumitem}
\usepackage{epstopdf}
\usepackage{algorithm}
\usepackage{algpseudocode}
\usepackage{array}
\usepackage{bbding}
\usepackage{fancyhdr} 
\usepackage{fancyvrb} %
\usepackage{longtable}
\usepackage{listings}
\usepackage{rotating,rotfloat} 
\usepackage{setspace}
\usepackage{subcaption}
\usepackage{tablefootnote}
\usepackage{threeparttable} 
\usepackage{float}
\usepackage{msc}
\usepackage[font=small,labelfont=bf,labelsep=none]{caption}

\newtheorem{remark}{Remark}[section]
\numberwithin{equation}{section} 

\begin{document}

\let\WriteBookmarks\relax
\def\floatpagepagefraction{1}
\def\textpagefraction{.001}
\shorttitle{}
\shorttitle{Robust and efficient solvers}
\shortauthors{Tan et al.}

\title [mode = title]{Robust and efficient solvers for nonlinear partial differential equations based on random feature method
}


\author[1,3]{Longze Tan}
\ead{ba23001003@mail.ustc.edu.cn}

\author[5]{Xiaohe Yue}
\ead{yuexiaohe666@gmail.com}

\author[1,2,3,4]{Jingrun Chen}
\ead{jingrunchen@ustc.edu.cn}

\author[2,3,4]{Jiamin Jiang}
\cormark[1]
\ead{jiaminjiang@ustc.edu.cn}

\affiliation[1]{School of Mathematical Sciences, University of Science and Technology of China, Hefei, China}

\affiliation[2]{School of Artificial Intelligence and Data Science, University of Science and Technology of China, Hefei, China}

\affiliation[3]{Suzhou Institute for Advanced Research, University of Science and Technology of China, Suzhou, China}

\affiliation[4]{Suzhou Big Data & AI Research and Engineering Center, Suzhou, China}

\affiliation[5]{School of Mathematical Sciences, East China Normal University,Shanghai, China}

\cortext[cor1]{Corresponding author}

\begin{abstract}
Solving nonlinear partial differential equations (PDEs) on complex geometries remains a fundamental yet challenging task in science and engineering. Traditional mesh-based methods often require carefully designed meshes and tailored approximation spaces, leading to high computational cost and implementation difficulties. 
The random feature method (RFM), a mesh-free machine
learning-based framework, has emerged as a promising alternative for solving PDEs on complex domains. However, for large three-dimensional nonlinear problems, attaining high accuracy typically requires domain partitioning with many collocation points and random features per subdomain, which leads to extremely large and ill-conditioned nonlinear least-squares systems. To overcome these challenges, we propose two randomized Newton-type solvers. The first is an inexact Newton method with right preconditioning (IPN), in which randomized Jacobian compression and QR factorization are used to construct an efficient preconditioner that substantially reduces the condition number. Each Newton step is then approximately solved by LSQR, and a derivative-free line search is incorporated to ensure residual reduction and stable convergence. Building upon this framework, we further develop an adaptive multi-step inexact preconditioned Newton method (AMIPN). In this approach, the preconditioned Jacobian is reused across multiple inner iterations, while a prescribed maximum number of inner iterations together with an adaptive early-stopping criterion determines whether the current preconditioner can be retained in subsequent outer iterations. These mechanisms effectively avoid redundant computations and enhance robustness. Extensive numerical experiments confirm the remarkable effectiveness of the proposed solvers. For three-dimensional steady-state PDEs, including elliptic, diffusion-reaction, Helmholtz, and Gray-Scott equations, the methods achieve substantially higher accuracy than classical finite element and finite difference approaches. For two-dimensional time-dependent problems with moving holes, topological changes, large deformations, and complex space-time geometries, the framework demonstrates robust convergence and notable improvements in accuracy. Compared with machine-learning-based solvers such as physics-informed neural networks and weak adversarial networks, our methods consistently yield order-of-magnitude gains, thereby establishing the RFM combined with IPN/AMIPN as an efficient framework for large-scale nonlinear PDEs.

\end{abstract}

\begin{keywords}
    Random feature method \sep  
    Nonlinear partial differential equations \sep
    Nonlinear solvers \sep 
    Adaptive\sep
    Preconditioning\sep
    Nonlinear Least Squares Problem

\end{keywords}

\maketitle
\section{Introduction}\label{sec:1}
The numerical solution of nonlinear partial differential equations (PDEs) plays a central role in numerous scientific and engineering applications. However, the problem becomes particularly challenging when the computational domain involves complex geometries, such as those undergoing topological changes, large deformations, or exhibiting intricate structures that render mesh generation infeasible. Classical numerical methods, including the finite difference method (FDM) \cite{leveque2007finite}, finite element method (FEM) \cite{zienkiewicz2005fem}, and spectral method \cite{shen2011spectral}, typically rely on carefully constructed approximation spaces and mesh structures tailored to the geometry, both of which require substantial computational effort and manual intervention  \cite{bueno2006spectralSISC, bueno2006spectral, li2009solving, yu2020higher, chen2017interface}. 
In recent years, deep learning-based approaches have emerged as promising alternatives for solving PDEs. Representative examples include physics-informed neural networks (PINNs) \cite{raissi2019physics}, the deep Galerkin method (DGM) \cite{sirignano2018dgm}, the weak adversarial networks (WAN) \cite{zang2020weak}, deep Ritz method \cite{EYu2018DeepRitz}, deep Nitsche method \cite{LiaoWang2021DeepNitsche}, and deep mixed residual method \cite{Lyu2022MIM}, as well as
other related approaches. The core idea is to approximate PDE solutions using deep neural networks, with physical laws and prior knowledge embedded in the loss function. These methods are inherently mesh-free, relatively simple to implement, and can flexibly handle irregular boundaries and evolving domains, thereby partially overcoming challenges arising from geometric complexity while also demonstrating notable advantages in tackling high-dimensional problems and inverse modeling tasks \cite{weinan2021dawning, karniadakis2021physics}.
\par
Despite their promising capabilities, deep neural network (DNN)-based methods often exhibit unsatisfactory performance due to the highly nonlinear and nonconvex nature of the underlying optimization problems. The training loss  are typically populated with numerous saddle points and sharp local minima, making gradient-based optimization particularly challenging. Common optimizers such as Adam \cite{adam2014method} and L-BFGS \cite{liu1989limited} frequently suffer from vanishing or exploding gradients and often stagnate at an $L^2$ relative error of about $1\%$ \cite{raissi2019physics, sirignano2018dgm, zang2020weak}, with little improvement even when the network depth or width is increased. As the model size increases, the difficulty of optimization also grows, which may lead to slower convergence and a risk of overfitting. For ill-conditioned or stiff PDEs, these challenges are more pronounced, as the loss are often steep \cite{wang2021understanding}, making it difficult for first-order methods to achieve high-accuracy solutions. The training process usually demands considerable time and memory resources and may deliver limited accuracy. Consequently, in many practical problems, current DNN-based PDE solvers have not yet fully reached the performance of traditional numerical methods such as the finite element or finite difference methods in terms of accuracy and efficiency. Nevertheless, DNN solvers demonstrate advantages for handling high-dimensional problems and complex geometries. Enhancing their convergence robustness, training efficiency, and accuracy will be an important direction toward making them practical alternatives in mainstream scientific computing.
\par
In low-dimensional settings, researchers have explored various strategies to integrate traditional numerical methods with machine learning techniques, aiming to leverage the strengths of both approaches \cite{yang2018novel, calabro2021extreme, dong2021local, CCEY}. A representative example is the random feature method (RFM) \cite{CCEY}, which constructs the global approximation space by combining partition of unity (PoU) with random feature functions. In this framework, random feature functions can be constructed as a two-layer neural network, where the inner-layer parameters are randomly generated and fixed, and only the output-layer parameters—the coefficients of the random features—are optimized. The PoU functions are typically formed via tensor products of univariate functions. For linear PDEs, the optimization problem associated with the RFM reduces to a convex least-squares system, whereas for nonlinear PDEs it becomes a nonlinear least-squares problem. Consequently, if these systems can be solved efficiently and accurately, the RFM is capable of delivering high-precision solutions without requiring full network training. This property offers an efficient alternative to conventional DNN training, particularly in scientific computing. The RFM has already been successfully applied to a broad spectrum of problems, including solid and fluid mechanics, interface and time-dependent problems, diffusion equations, elliptic eigenvalue problems, and topology optimization \cite{CCEY, CCYang, CELuo, mei2024solving, xiong2025random, mei2025direct}. Moreover, the framework has been extended to address more challenging scenarios: for instance, a two-level RFM has been developed by incorporating domain decomposition for elliptic PDEs on geometrically complex domains \cite{sun2025two}, and an asymptotic-preserving RFM has been proposed by integrating micro-macro decomposition strategies for multiscale radiative transfer equations \cite{chen2025micro}.
\par
In this work, we first introduce the basic framework of the RFM for solving partial differential equations (PDEs). We then observe that when classical nonlinear solvers—such as Newton's method~\cite{kelley1999iterative}, the Levenberg-Marquardt (LM) method~\cite{levenberg1944method,marquardt1963algorithm}, Mor{\'e}'s algorithm~\cite{more1978lm} (as implemented in the MINPACK package), and trust-region (TR) methods~\cite{yuan1994trust}—are applied to the large-scale nonlinear least-squares (NLS) problems arising from RFM discretizations of three-dimensional nonlinear PDEs, the main difficulties stem from the fact that the associated Jacobian matrices are often highly ill-conditioned. As a consequence, the linear systems required in the nonlinear iterations can be difficult to solve efficiently and stably at large scales, which may significantly affect the overall convergence performance. 
\par
To enable the practical application of the RFM to large-scale nonlinear PDEs, this study proposes two efficient and robust randomized Newton-type solvers. we first propose an inexact Newton method with right preconditioning, where at each iteration the Jacobian is compressed via count sketch, followed by QR factorization to construct a right preconditioner. The preconditioned system is then solved by LSQR, combined with a derivative-free line search to guarantee residual reduction and convergence stability. Secondly, to further accelerate convergence, building upon IPN, we further propose an adaptive multi-step inexact preconditioned Newton method. In this framework, each outer iteration requires only one Jacobian evaluation and preconditioning, and the resulting preconditioned Jacobian can be reused across multiple inner Newton iterations. To avoid redundant computations, the inner loop is assigned a fixed maximum number of inner iterations and equipped with an adaptive early-stopping criterion to determine whether the current preconditioned Jacobian can be reused in the next outer iteration. Moreover, unlike IPN, AMIPN incorporates a derivative-free line search within the inner loop to ensure residual decrease at every update, thereby further enhancing convergence stability.

Main contributions of this work are summarized as follows:
\begin{itemize}
\item We establish the combination of RFM with the proposed IPN/AMIPN solvers as an efficient and robust mesh-free computational framework for solving nonlinear PDEs on complex geometries. The proposed nonlinear solvers can effectively handle overdetermined and ill-conditioned large-scale systems arising in scientific and engineering applications, and provide a potential pathway for the future numerical simulation of complex physical processes and coupled problems. 
\item Comprehensive experiments demonstrate that the solvers achieve accuracy improvements of several orders of magnitude compared with finite element, finite difference, and learning-based approaches (PINNs and WANs), while maintaining strong robustness and rapid convergence across a broad spectrum of steady-state and time-dependent PDEs with complex and evolving geometries.
\item To the best of our knowledge, this is the first application of randomized preconditioning techniques to overdetermined large-scale nonlinear systems from neural-network-based discretizations of nonlinear PDEs, thereby constructing efficient and robust nonlinear solvers.
\end{itemize}
\par
The remainder of this paper is organized as follows. Section~\ref{sec:2} introduces the class of nonlinear PDEs under consideration and provides a brief overview of the PINNs and RFM. The convergence limitations of classical nonlinear solvers are illustrated in Section~\ref{subsec:3.4} through a representative example. Sections~\ref{sec:4} and \ref{sec:5} form the core of the paper: the former presents two new nonlinear solvers, IPN and AMIPN, whereas the latter develops the overall RFM-based framework that incorporates these solvers for the efficient solution of nonlinear PDEs. Section \ref{sec:6} reports numerical results on three-dimensional steady-state and two-dimensional time-dependent problems with complex geometries. Section \ref{sec:7} concludes the paper, and Appendix \ref{Appendix_A} provides additional nonlinear examples.

\section{Preliminaries}\label{sec:2}
In this section, we first present the mathematical formulation of the nonlinear PDE under consideration, followed by the introduction of the strong-form-based PINNs and RFM.
\subsection{The nonlinear PDE problem}
We consider the following general nonlinear partial differential equation: find $u: \Omega \rightarrow \mathbb{R}$ such that
\begin{equation}\label{boundary-value problem}
\left\{
\begin{aligned}
\mathcal{P}(u(\mathbf{x})) &= f(\mathbf{x}), & \mathbf{x} \in \Omega, \\
\mathcal{B} u(\mathbf{x}) & =g(\mathbf{x}), & \mathbf{x} \in \partial \Omega.
\end{aligned}
\right.
\end{equation}
where $\mathcal{P}(u(\mathbf{x}))=\lambda \mathcal{L} u(\mathbf{x})+\mu \mathcal{N}(u(\mathbf{x}))$. Here, $\mathcal{L}$ and $\mathcal{N}$ denote a linear and a nonlinear differential operator, respectively; $u(\mathbf{x}) \in \mathbb{R}$ is the unknown field; $\mathcal{B}$ is a (possibly differential or algebraic) boundary operator; $f$ and $g$ are prescribed functions; and $\lambda, \mu \in \mathbb{R}$ are constants. When $\mu=0$, the problem reduces to a linear one; throughout this work we assume $\mu \neq 0$.
\par
For time-dependent problems we allow $\mathcal{L}$ to involve temporal derivatives (e.g., $\partial_t, \partial_{t t}$ ). We treat time on the same footing as space by augmenting the spatial coordinates with the temporal variable, and regard the problem as $d$-dimensional with $\mathbf{x}=\left(x_1, \ldots, x_{d-1}, x_d=t\right) \in \Omega \subset \mathbb{R}^{d-1} \times[0, T]$ Accordingly, the computational domain $\Omega$ encompasses both spatial and temporal dimensions, and the boundary operator $\mathcal{B}$ encodes not only the spatial boundary conditions on $\partial \Omega \times(0, T]$ but also the initial condition on $\Omega \times\{0\}$, the latter being interpreted as a special type of boundary condition in the spatiotemporal domain.

\subsection{Physics-informed neural networks}
To solve problem \eqref{boundary-value problem}, the physics-informed neural network is employed to approximate the target solution $u(\mathbf{x})$ :
\begin{align*}
    u(\mathbf{x}) \approx \hat{u}(\mathbf{x} ; \Theta) \triangleq \mathcal{N}(\mathbf{x} ; \Theta)
\end{align*}
where $\mathcal{N}(\mathbf{x}; \Theta)$ denotes a feed-forward neural network with one input layer, $L-1$ hidden layers, and one output layer, defined as
\begin{align*}
    \begin{cases}\mathbf{z}^1=W^1\mathbf{x}+\mathbf{b}^1, & \mathbf{h}^1=\sigma\left(\mathbf{z}^1\right), \\ \mathbf{z}^l=\mathbf{W}^l \mathbf{h}^{l-1}+\mathbf{b}^l, &\mathbf{h}^l=\sigma\left(\mathbf{z}^l\right), \quad l=2,3, \ldots, L-1 \\ \mathcal{N}(\mathbf{x}; \Theta)=\mathbf{W}^L \mathbf{h}^{L-1}+\mathbf{b}^L, & \end{cases}
\end{align*}
with $\sigma(\cdot)$ being the nonlinear activation function, and $\mathbf{W}^l, \mathbf{b}^l$ representing the weight matrix and bias vector between the $(l-1)$-th and $l$-th layers, respectively. The set of trainable parameters is given by
\begin{align*}
    \Theta=\bigcup_{l=1}^L\left\{\mathbf{W}^l, \mathbf{b}^l\right\}
\end{align*}
which collects all weights and biases in the network.

The network parameters $\Theta$ are obtained by minimizing a physics-informed loss function consistent with the governing equation, defined as
\begin{align*}
    \mathcal{L}_{\text {total }}(\hat{u} ; \Theta)=\mathcal{L}_{\text {Data }}+\mathcal{L}_{\mathrm{PDE}}+\mathcal{L}_{\mathrm{BC}}
\end{align*}
where
\begin{align*}
\mathcal{L}_{\text {Data }}&=\lambda_d \sum_{i=1}^{N_d}\left|\hat{u}\left(\mathbf{x}_i ; \Theta\right)-u_i\right|^2, \quad \mathbf{x}_i \in \hat{D} \\
\mathcal{L}_{\mathrm{PDE}}&=\lambda_r \sum_{i=1}^{N_{ \hat{\Omega}}}\left|\mathcal{P}\left(\hat{u}\left(\mathbf{x}_i ; \Theta\right)\right)-f\left(\mathbf{x}_i\right)\right|^2=\lambda_r \sum_{i=1}^{N_{\hat{\Omega}}}\left|\lambda \mathcal{L} \hat{u}\left(\mathbf{x}_i ; \Theta\right)+\mu \mathcal{N}\left(\hat{u}\left(\mathbf{x}_i ; \Theta\right)\right)-f\left(\mathbf{x}_i\right)\right|^2, \quad \mathbf{x}_i \in  \hat{\Omega} \\
\mathcal{L}_{\mathrm{BC}}&=\lambda_b \sum_{i=1}^{N_{\partial  \hat{\Omega}}}\left|\mathcal{B}\left(\hat{u}\left(\mathbf{x}_i ; \Theta\right)\right)-g\left(\mathbf{x}_i\right)\right|^2, \quad \mathbf{x}_i \in \partial  \hat{\Omega}
\end{align*}
Here, $\hat{D}$ is the set of labeled samples of size $N_d ; \hat{\Omega} \subset \Omega$ is the set of collocation points with size $N_{\hat{\Omega}}$; and $\partial \hat{\Omega} \subset \partial \Omega$ is the set of boundary points with size $N_{\partial \hat{\Omega}}$. The coefficients $\lambda_d, \lambda_r$, and $\lambda_b$ are nonnegative weighting factors that balance the contributions of the data, PDE residual, and boundary condition terms, respectively. The differential operators $\mathcal{L}$ and $\mathcal{B}$ are evaluated by means of automatic differentiation. Although PINNs provide a promising framework for solving PDEs, several challenges remain to be addressed: (i) potential imbalances in the multi-objective loss due to scale differences among data, PDE residual, and boundary terms; (ii) possible limitations in accuracy and efficiency resulting from the inherent non-convex optimization, which may lead to convergence to suboptimal solutions; and (iii) a certain reliance on carefully designed architectures or adaptive strategies, which can increase implementation complexity.

\subsection{The random feature method}\label{sec:2.2}
To overcome the limitations of PINNs, single-hidden-layer neural network approaches such as LocELM \cite{dong2021local} and RFM have been proposed, in which the hidden-layer parameters are fixed and only the output weights are trained. LocELM, however, often encounters imbalance in problems such as structural mechanics, where the magnitudes of the governing equations and boundary conditions differ substantially, leading to limited accuracy \cite{CCEY}. In contrast, RFM incorporates a rescaling strategy in the loss function to balance the contributions of PDE residuals and boundary/initial conditions, thereby improving both convergence and accuracy.

In RFM, the approximate solution $u_n(\mathbf{x})$ is expressed as a linear combination of $n$ random feature functions defined on $\Omega$, i.e., 
\begin{equation}\label{RFM_approximation_solution}
  u_n(\mathbf{x})=\sum_{j=1}^n u_j \phi_j(\mathbf{x}).
\end{equation}
While random feature functions are globally defined, solutions to PDEs often exhibit localized variations, potentially at small scales. To capture such behavior, we construct multiple local approximations, each based on an independent random feature model, and blend them seamlessly using partition of unity (PoU) functions.
\par
To construct the PoU, we first introduce a hypercube
$
\Omega_c=\prod_{i=1}^d\left[a_i, b_i\right]
$
of appropriate size that fully contains the domain $\Omega$. Here, $\prod_{i=1}^d\left[a_i, b_i\right]$ denotes the Cartesian product of $d$ intervals, i.e., $\left[a_1, b_1\right] \times\left[a_2, b_2\right] \times \cdots \times\left[a_d, b_d\right]$.
The hypercube $\Omega_c$ is then partitioned into $M_p$ non-overlapping hyperrectangles $\left\{\Omega_i\right\}_{i=1}^{M_p}$ as
$
    \Omega_c=\bigcup_{i=1}^{M_p} \Omega_i,~ \Omega_i=\prod_{j=1}^d\left[a_{i j}, b_{i j}\right],
$
with adjacent elements satisfying $a_{i+1, j}=b_{i j}$. This construction is analogous to structured meshing in the finite element method, and can be generalized to more complex partitioning schemes if needed.
\par
Let $\mathbf{x}=\left(x_1, \ldots, x_d\right)^T \in \Omega$. For each subdomain $\Omega_i$, we introduce the following affine transformation to obtain a partition-dependent normalized coordinate vector $\mathbf{l}_i(\mathbf{x})=\left(\tilde{x}_{i 1}, \ldots, \tilde{x}_{i d}\right)^T$ :
\begin{align*}
    \tilde{x}_{i j}=2 \frac{x_j-a_{i j}}{b_{i j}-a_{i j}}-1, \quad i=1, \ldots, M_p, \quad j=1, \ldots, d.
\end{align*}
Denote the center and half-width (radius) of the hyperrectangle $\Omega_i$ as
\begin{align*}
    \boldsymbol{\mu}_i=\left(\frac{b_{i 1}+a_{i 1}}{2}, \ldots, \frac{b_{i d}+a_{i d}}{2}\right)^T, \quad \boldsymbol{\sigma}_i=\left(\frac{b_{i 1}-a_{i 1}}{2}, \ldots, \frac{b_{i d}-a_{i d}}{2}\right)^T.
\end{align*}
Using these notations, the transformation can be compactly expressed in vector form as
\begin{align*}
    \mathbf{l}_i(\mathbf{x})=\frac{\mathbf{x}-\boldsymbol{\mu}_i}{\boldsymbol{\sigma}_i}, \quad i=1, \ldots, M_p.
\end{align*}
This transformation is a simple coordinate scaling and translation that maps the physical subdomain $\Omega_i$ bijectively onto the reference hypercube $[-1,1]^d$. Thus, the PoU functions $\left\{\psi_i(\mathbf{x})\right\}_{i=1}^M$ are constructed via the tensor product of a univariate function $\varphi$ :
\begin{align*}
    \psi(\mathbf{x})=\prod_{j=1}^d \varphi\left(\tilde{x}_{i j}\right), \quad i=1, \ldots, M_p
\end{align*}
where $\varphi$ is selected following \cite{CCEY} as
\begin{equation}
\varphi^a(y) =
\begin{cases}
1, & |y| \leq 1, \\
0, & \text{else},
\end{cases}
\quad\text{or}\quad
\varphi^b(y) =
\begin{cases}
1, & |y| \leq \frac{3}{4}, \\[2mm]
\dfrac{1 - \sin(2\pi |y|)}{2}, & \frac{3}{4} \leq |y| \leq \frac{5}{4}, \\[2mm]
0, & \text{else}.
\end{cases}
\end{equation}
$\varphi_a(y)$ is discontinuous with a smaller support, while $\varphi_b(y)$ is continuously differentiable with a larger support. 
\par
Next, $J_i$ random feature functions on each subdomain $\Omega_i$ are constructed as follows
\begin{align}
  \phi_{ij}(\mathbf{x})=\sigma\left(\mathbf{k}_{ij} \cdot \mathbf{l}_i(\mathbf{x})+b_{ij}\right), \quad i=1, \ldots, J_i,
\end{align}
where the nonlinear activation function $\sigma$ is often chosen as tanh or trigonometric functions and each component of $\mathbf{k}_{ij}$ and $b_{ij}$ is chosen uniformly from the interval $\left[-R_{ij}, R_{ij}\right]$ and is fixed. The parameter $R_{ij}$ is used to control the initialization range of the weights $\left\{\mathbf{k}_{ij} \right\}$ and $\left\{b_{i j}\right\}$. Thus, the approximate solution $u_n(\mathbf{x})$ can be represented as 
\begin{equation}\label{RFM_approximate_sol}
    u_n(\mathbf{x})=\sum_{i=1}^{M_p} \psi_i(\mathbf{x}) \sum_{j=1}^{J_i} u_{ij} \phi_{ij}(\mathbf{x}),
\end{equation}
where $u_{ij}\;(i=1,\ldots,M_p, j=1,\ldots, J_i)$ are the unknown coefficients that need to be determined, and $n = (J_1 + \ldots+ J_{M_p})$ denotes the degrees of freedom. 
A loss function for \eqref{boundary-value problem} can be formulated as follows
\begin{equation}\label{loss_function}
  Loss=\sum_{\mathbf{x}_i \in C_I} \sum_{k=1}^{K_I} \lambda_{I i}^k\|\mathcal{P}^k (u_n\left(\mathbf{x}_i\right))-f^k\left(\mathbf{x}_i\right)\|_{l^2}^2+\sum_{\mathbf{x}_j \in C_B} \sum_{\ell=1}^{K_B} \lambda_{B j}^{\ell}\|\mathcal{B}^{\ell} u_n\left(\mathbf{x}_j\right)-g^{\ell}\left(\mathbf{x}_j\right)\|_{l^2}^2. 
\end{equation}
Here $ \lambda_{I i}^k$ and $\lambda_{B j}^{\ell}$ are penalty parameters, as defined in \cite{CCEY} and $\left\{\mathbf{x}_i\right\} \subset \Omega$ and $\left\{\mathbf{x}_j\right\} \subset \partial \Omega$ are collocation points. 
It is worth noting that, when $M_p>1$ and PoU function $\varphi^a(y)$ is used, smoothness conditions between the adjacent elements in the partition are explicitly imposed by adding regularization terms in loss function \eqref{loss_function} while no regularization is required when $\varphi^b(y)$ is used for second-order equations due to its first-order continuity. Next, we let $K_I$ and $K_B$ be the number of condition at each interior point and boundary point, respectively. The total number of conditions is $m=K_I \# C_I + K_B \# C_B$. 
In practice, the coefficients $u_{i j}$ are the only unknown variables in \eqref{loss_function} when $\mu=0$. Therefore, the optimization problem \eqref{loss_function} can be regarded as a NLS problem with respect to $u_{i j}$.

\begin{remark}
  Within the RFM framework, for vector-valued solutions, we approximate each component of the solution individually using \eqref{RFM_approximate_sol}, i.e.,
\begin{align*}
    \boldsymbol{u}_n(\mathbf{x})=\left(u_n^1(\mathbf{x}), \ldots, u_n^{K_I}(\mathbf{x})\right)^{\top},
\end{align*}
where $u_n^q(\mathbf{x})=\sum_{i=1}^{M_p} \psi_i(\mathbf{x}) \sum_{j=1}^{J_i} u_{i j}^q \phi_{i j}(\mathbf{x}),~q=1, \ldots, K_I$, and $K_I$ denotes the output dimension.
\end{remark}

\section{Classical nonlinear solvers}\label{subsec:3.4}
In this section, we first review several classical nonlinear solvers and then illustrate their numerical results through a three-dimensional nonlinear elliptic PDE example, thereby providing motivation for the methods proposed in this work.
\subsection{Nonlinear solvers}
Newton's method and the Gauss-Newton method are the most fundamental  approaches for solving \eqref{loss_function}, and many efficient NLS algorithms have been developed based on them. Since \eqref{loss_function} is overdetermined (the number of equations exceeds the number of unknowns), the Newton update at the $(k+1)$-th iteration is given by
\begin{equation}\label{newton_iterate}
    \mathbf{u}_{k+1}=\mathbf{u}_k-\mathbf{J}_k^{\dagger} \mathbf{F}_k,
\end{equation}
where $\mathbf{u}_k$ denotes the $k$-th iterate, $\mathbf{F}_k=\mathbf{F}(\mathbf{u}_k)$, $\mathbf{J}_k=\mathbf{F}^{\prime}(\mathbf{u}_k)$ is the Jacobian evaluated at $\mathbf{u}_k$, and $\mathbf{J}_k^{\dagger}$ is the Moore-Penrose pseudoinverse of $\mathbf{J}_k$. The iteration \eqref{newton_iterate} can be derived from the following linear least-squares problem:
\begin{equation}\label{newton_system}
\min_{\mathbf{u} \in \mathbb{R}^n}\|\mathbf{F}_k+\mathbf{J}_k(\mathbf{u}-\mathbf{u}_k)\|_2,
\end{equation}
which corresponds to minimizing the residual norm of the first-order Taylor expansion of $\mathbf{F}(\mathbf{u})$ at $\mathbf{u}_k$. By solving the normal equations corresponding to \eqref{newton_system}, one obtains the Gauss-Newton (GN) method.

The GN method requires the Jacobian $\mathbf{J}_k$ at $\mathbf{u}_k$ to have full column rank; otherwise the update may be undefined. When $\mathbf{J}_k$ is ill-conditioned, the GN step can be excessively long and cause instabilities. To remedy this, Levenberg \cite{levenberg1944method} and Marquardt \cite{marquardt1963algorithm} introduced a nonnegative damping parameter $\lambda_k$, yielding the iteration:
\[
    \mathbf{u}_{k+1}=\mathbf{u}_k-(\mathbf{J}_k^{\top} \mathbf{J}_k+\lambda_k \mathbf{I})^{-1} \mathbf{J}_k^{\top} \mathbf{F}_k .
\]
The Levenberg-Marquardt (LM) method is always well defined when $\lambda_k>0$. Common choices include $\lambda_k=\|\mathbf{F}_k\|^2$ \cite{yamashita2001rate}, which ensures local quadratic convergence but may suffer from numerical difficulties, and $\lambda_k=\|\mathbf{F}_k\|$ \cite{fan2005quadratic}, which alleviates these issues.
\par
Beyond LM-type damping, a representative example is the well-known Mor{\'e} algorithm \cite{more1978lm}, which is implemented in the MINPACK software package. Moreover, The LM method can be interpreted as a special case of trust-region approaches, where the damping parameter implicitly controls the step length instead of updating the trust-region radius explicitly \cite{yuan1994trust,more1983computing}.

\subsection{Difficulties of classical nonlinear solvers}
In this section, we consider a benchmark example to present the numerical results of classical nonlinear solvers, including Newton, LM, Mor{\'e}, and TR methods, when applied to \eqref{loss_function}. Specifically, we examine the following three-dimensional  nonlinear elliptic problem on the unit cube $\Omega=(0,1)^3$:
\begin{equation}\label{benchmark_nonlinear_elliptic_equation}
\begin{cases}-\nabla \cdot(\Lambda(u) \nabla u)=f(x, y, z), & (x, y, z) \in \Omega  \\  \Lambda(u)=1+u^2 & \end{cases}
\end{equation}
The exact solution is given by
\begin{align*}
     u(x, y, z)=(x+y+z) \sin (\pi x) \sin (\pi y) \cos (\pi z)
\end{align*}
and the corresponding Dirichlet boundary conditions are imposed on $\partial\Omega$.
\par
We uniformly partition $\Omega$ along each coordinate direction into $N_x$, $N_y$, and $N_z$ subdomains, resulting in $M_p=N_xN_yN_z$ subdomains. Each subdomain is denoted by $\Omega_{e_{mns}}$with indices $0\le m\le N_x-1$, $0\le n\le N_y-1$, and $0\le s\le N_z-1$, where $e_{mns}$ serves as the subdomain label. On each $\Omega_{e_{mns}}$, collocation points are placed as $(x_p^{e_{mns}},y_q^{e_{mns}},z_r^{e_{mns}})$ with $0\le p\le Q_x-1$, $0\le q\le Q_y-1$, and $0\le r\le Q_z-1$, where ${x_p^{e_{mns}}}$, ${y_q^{e_{mns}}}$, and ${z_r^{e_{mns}}}$ denote sets of $Q_x$, $Q_y$, and $Q_z$ points on the respective intervals. With $Q=Q_xQ_yQ_z$ points per subdomain, the total number of interior collocation points is $Q_I=N_xN_yN_zQ$. For problem \eqref{benchmark_nonlinear_elliptic_equation}, we adopt uniformly spaced collocation points including endpoints. Each subdomain is assigned the same number $J$ of random feature functions. In the numerical experiments, $\Omega$ is partitioned into $M_p=8$ subdomains with $N_x=N_y=N_z=2$, collocation points are set with $Q_x=Q_y=Q_z\in \left\{20,25,30\right\}$, and the number of random features $J$ is varied as $800$, $1200$, $1600$, and $2000$. The PoU function $\psi_i^a(\mathbf{x})$ is employed, the activation function is chosen as \texttt{tanh}, and each component of $\mathbf{k}_{ij}$ and $b_{ij}$ is independently sampled from the uniform distribution on $[-1,1]$ and fixed throughout the computations. Under these settings, solving \eqref{benchmark_nonlinear_elliptic_equation} with the RFM reduces to the NLS problem
\begin{equation}\label{NLS_problem}
\min_{\mathbf{u}\in\mathbb{R}^n}\|\mathbf{F}(\mathbf{u})\|_2,
\end{equation}
where $\mathbf{F}:\mathbb{R}^n\to\mathbb{R}^m$, with $n=8J$ unknown coefficients and $m=Q_I+Q_b+Q_c$. Here, $Q_b=2(N_xN_yQ_xQ_y+N_xN_zQ_xQ_z+N_yN_zQ_yQ_z)$ corresponds to the number of Dirichlet boundary constraints, and $Q_c=2\big((N_x-1)N_yN_zQ_yQ_z+(N_y-1)N_xN_zQ_xQ_z+(N_z-1)N_xN_yQ_xQ_y\big)$ arises from the regularization enforcing smoothness across adjacent subdomains.
\par
All numerical experiments were performed on a computing server equipped with two Intel Xeon Gold 6342 CPUs (48 cores, 96 threads, 2.80 GHz base and 3.50 GHz turbo), 512 GB of RAM, and four NVIDIA A800 GPUs (80 GB memory each, PCIe interface). 
The implementation was developed in Python 3.10.18 with CUDA 12.2 support, and executed under a 64-bit Linux operating system.

Within this computational environment, we investigate four nonlinear solvers for the NLS problem \eqref{NLS_problem}: the Newton, LM, Mor{\'e}, and TR methods. The latter two are implemented in \texttt{scipy.optimize.least\_squares}. The numerical results are summarized in Table~\ref{3D_nonlinear_equation_newton_LM}, where ``IT'' denotes the number of iterations and ``CPU'' refers to the elapsed CPU time for solving \eqref{NLS_problem}. The stopping criterion is chosen as $\|\mathbf{J}_k^T \mathbf{F}_k\|<10^{-6}$, with maximum iterations set to 60 for Newton and 180 for LM. In the implementation, Newton updates the solution using \texttt{scipy.linalg.lstsq}, which calls the LAPACK routine \texttt{gelsd} based on a divide-and-conquer SVD algorithm, whereas LM employs a Cholesky decomposition with the damping parameter chosen as $\lambda_k=\|\mathbf{F}_k\|$ at each iteration. For PDE problems, accuracy is evaluated using the relative $L^2$ error with respect to the exact solution.
\par
Overall, due to the highly ill-conditioned and large-scale nature of the Jacobian matrices, classical solvers fail to solve the Newton linear systems both efficiently and accurately, leading to pronounced convergence difficulties in this example. In particular, the Mor{\'e} method did not converge within the 12-hour time limit even for the smallest problem size. As the number of collocation points and basis functions per subdomain increased, the computational costs of Newton, LM, and TR methods grew drastically. Moreover, Newton and LM often failed to converge within the prescribed maximum number of iterations, with their $u$-errors stagnating at approximately $10^{-4}$ and $10^{-2}$, respectively. By contrast, the TR method achieved higher accuracy under certain configurations. For example, when $\sqrt[3]{Q}=25$ with $J=1200$ and $1600$, the relative $L^2$ errors were $1.32 \times 10^{-6}$ and $1.18 \times 10^{-7}$, respectively—significantly better than the $10^{-4}$-level errors obtained by Newton. However, these accuracy gains came at the expense of substantially increased computational cost. For larger-scale systems, the TR method further failed in the underlying SVD routine due to excessive memory requirements.

\begin{table}[!htbp]
  \centering
  \begin{threeparttable}
  \caption{Comparison of the Newton and LM methods for the 3D nonlinear elliptic problem \eqref{benchmark_nonlinear_elliptic_equation}}
  \label{3D_nonlinear_equation_newton_LM}
  \begin{tabular}{c c c c c c c c}
    \toprule
    $(N_x, N_y, N_z)$ & $\sqrt[3]{Q}$ & $J$ & Solver & IT & CPU & $u$ error & $H^1$ error \\
    \midrule
    \multirow{15}{*}{$(2, 2, 2)$} & \multirow{4}{*}{$20$} 
        & \multirow{4}{*}{800} & Newton & $11$ & $1280.39$  & $3.13e-05$ & $5.11e-04$ \\
        & & & LM     & $180$ & $ 2656.61$  & $1.00e-01$ & $6.79e-01$ \\
        & & & TR     & $57$ & $6607.49$  & $3.13e-05$ & $5.11e-04$ \\
        & & & Mor{\'e}      & \multicolumn{4}{c}{\textemdash\hyperlink{note:a}{\textsuperscript{a}}} \\
    \cmidrule(lr){2-8}
        & \multirow{6}{*}{$25$} 
        & \multirow{3}{*}{1200} & Newton & $60$ & $7076.28$  & $2.45e-04$ & $5.41e-03$ \\
        & & & LM     & $180$ & $7139.62$  & $ 6.14e-02$ & $4.59e-01$ \\
         & & & TR     & 55 & $22745.37$ & $1.32e-06$ & $2.97e-05$ \\
        \cmidrule(lr){4-8}
        && \multirow{3}{*}{1600} & Newton & $60$ & $12484.37$  & $2.75e-04$ & $6.64e-03$ \\
        & & & LM     & $180$ & $13768.50$  & $5.61e-02$ & $4.28e-01$ \\
        & & & TR     & $57$ & $46132.25$ & $1.18e-07$ & $3.19e-06$ \\
      \cmidrule(lr){2-8}
        & \multirow{5}{*}{$30$} 
        & \multirow{3}{*}{1600} & Newton & $60$ & $18426.67$  & $3.37e-04$ & $7.72e-03$ \\
        & & & LM     & $180$ & $19728.05$  & $5.01e-02$ & $3.88e-01$ \\
        & & & TR     & \multicolumn{4}{c}{\textemdash\hyperlink{note:b}{\textsuperscript{b}}} \\
        \cmidrule(lr){4-8}
        && \multirow{2}{*}{2000} & Newton & $60$ & $ 28162.52$ & $3.53e-04$ & $8.64e-03$ \\
        & & & LM     & \multicolumn{4}{c}{\textemdash\hyperlink{note:c}{\textsuperscript{c}}} \\
    \bottomrule
  \end{tabular}

  \begin{tablenotes}[flushleft]
  \footnotesize
  \item[\hypertarget{note:a}{}\hyperlink{note:a}{\textsuperscript{a}}]
    The Mor{\'e} method in \texttt{scipy.optimize.least\_squares} did not converge within the 12-hour time limit, and the computation was therefore terminated.

  \item[\hypertarget{note:b}{}\hyperlink{note:b}{\textsuperscript{b}}]
  A numerical breakdown was observed when employing the TR method within the \texttt{scipy.optimize.least\_squares} routine. Specifically, the underlying LAPACK implementation failed during the SVD step due to excessive memory requirements. The routine reported the following runtime exception:
  \texttt{ValueError: Indexing a matrix of 3317760000 elements would incur an integer overflow in LAPACK.}

  \item[\hypertarget{note:c}{}\hyperlink{note:c}{\textsuperscript{c}}]
    The runtime of the LM method exceeds 24 hours, and the iteration is terminated.
  \end{tablenotes}

  \end{threeparttable}
\end{table}

\section{Proposed nonlinear solvers}\label{sec:4}
The results in the previous section indicate that solving the large-scale and highly ill-conditioned NLS problems arising from RFM discretizations of nonlinear PDEs requires the development of efficient and robust solvers. To this end, we propose two randomized Newton-type methods, namely an inexact Newton method with right preconditioning and its adaptive multi-step extension.
\subsection{The IPN method}

The underlying optimization problem \eqref{loss_function} can be written as
\begin{equation}\label{RFM_NLS}
\min_{\mathbf{u}\in \mathbb{R}^n}\|\mathbf{F}(\mathbf{u})\|_2,
\end{equation}
where the $k$-th Newton correction $\delta \mathbf{u}_k$ is obtained from the linearized system
\begin{equation}\label{RFM_Newton_system}
\mathbf{J}(\mathbf{u}_k)\delta \mathbf{u}=-\mathbf{F}(\mathbf{u}_k),
\end{equation}
with $\mathbf{J}(\mathbf{u}_k)\in \mathbb{R}^{m\times n}$ denoting the Jacobian. For large-scale and ill-conditioned systems, solving \eqref{RFM_Newton_system} dominates the overall cost and often deteriorates the convergence of Newton-type methods.

To stabilize the iteration, we introduce a right preconditioner $\mathbf{R}_k^{-1}$ and reformulate the Newton system as
\begin{equation}
\mathbf{J}(\mathbf{u}_k)\mathbf{R}_k^{-1}\mathbf{y}=-\mathbf{F}(\mathbf{u}_k),
\qquad \delta \mathbf{u}=\mathbf{R}_k^{-1}\mathbf{y}.
\end{equation}
The preconditioner $\mathbf{R}_k$ is constructed such that the preconditioned Jacobian $\tilde{\mathbf{J}}_k=\mathbf{J}(\mathbf{u}_k)\mathbf{R}_k^{-1}$ exhibits a substantially reduced condition number with well-clustered singular values. This ensures fast and stable convergence of Krylov subspace solvers for the inner least-squares problems, thereby enabling efficient and robust Newton iterations.

We employ a randomized sketching strategy to construct $\mathbf{R}_k$: (i) the Jacobian is compressed via count sketch \cite{ChChFC,ClWo,Woodruff},
\begin{align*}
    \tilde{\mathbf{B}}_k=\mathbf{S}_k \mathbf{J}\left(\mathbf{u}_k\right), \quad \mathbf{S}_k \in \mathbb{R}^{s \times m}, s=\gamma n, \gamma>1 ;
\end{align*}
(ii) a thin QR factorization is computed as
\begin{align*}
    \tilde{\mathbf{B}}_k=\mathbf{Q}_k \mathbf{R}_k, \quad \mathbf{R}_k \in \mathbb{R}^{n \times n}
\end{align*}
The computation of $\tilde{\mathbf{B}}_k$ requires only $\mathcal{O}(\mathrm{nnz}(\mathbf{J}_k))$ operations, and in practice the count sketch matrix $\mathbf{S}_k$ is never formed explicitly, thereby substantially reducing both memory usage and computational cost.

In summary, the first baseline nonlinear solver, referred to as the IPN method, is outlined in Algorithm~\ref{alg:irpn}.
\begin{algorithm}[!htbp]
\caption{IPN method}
\label{alg:irpn}
\begin{algorithmic}[1]
\Require Nonlinear mapping $\mathbf{F}(\mathbf{u})$, Jacobian $\mathbf{J}(\mathbf{u})$, initial guess $\mathbf{u}_0$, oversampling factor $\gamma>1$, tolerances $\varepsilon,\eta$
\Ensure approximate solution $\mathbf{u}$
\For{$k=0,1,2,\dots$ until convergence}
    \State The preconditioner $\mathbf{R}_k$ is obtained by Count Sketch compression of $J(\mathbf{u}_k)$ to $\tilde{\mathbf{B}}_k$ followed by thin QR factorization.
    \State Form the preconditioned matrix $\tilde{\mathbf{J}}_k=\mathbf{J}\left(\mathbf{u}_k\right) \mathbf{R}_k^{-1}$.
    \State Solve the preconditioned system $\min _{\mathbf{y} \in \mathbb{R}^n}\|\tilde{\mathbf{J}}_k \mathbf{y}+\mathbf{F}\left(\mathbf{u}_k\right)\|_2$ approximately using LSQR with tolerance $\eta$.
    \State Compute the Newton direction $\delta \mathbf{u}_k=\mathbf{R}_k^{-1}\mathbf{y}_k$
    \State Choose $\alpha_k>0$ by golden-section search such that $\|\mathbf{F}(\mathbf{u}_k+\alpha_k \delta \mathbf{u}_k)\|$ sufficiently decreases.
    \State Set $\mathbf{u}_{k+1}=\mathbf{u}_k+\alpha_k \delta \mathbf{u}_k$.
    \If{$|\|\mathbf{F}(\mathbf{u}_{k+1})\| - \|\mathbf{F}(\mathbf{u}_k)\|| < \varepsilon$}
    \State \textbf{break}
  \EndIf
\EndFor
\State \Return $\mathbf{u} \gets \mathbf{u}_{k+1}$
\end{algorithmic}
\end{algorithm}

It was shown in \cite{chen2024high} that, under mild assumptions, the preconditioned system satisfies
\begin{align*}
\kappa\left(\tilde{\mathbf{J}}_k\right) \leq \sqrt{\frac{1+\epsilon}{1-\epsilon}}, \qquad \epsilon,\delta \in (0,1),
\end{align*}
with probability at least $1-\delta$. Importantly, this bound is independent of the condition number $\kappa\left(\mathbf{J}(\mathbf{u}_k)\right)$. Numerical experiments further confirm that the proposed approach significantly improves the singular value distribution, effectively transforming an ill-conditioned system into a well-conditioned one. As a result, the right-preconditioning strategy provides stable and efficient IPN iterations for large-scale NLS problems.

\subsection{The AMINP method}
In Algorithm~\ref{alg:irpn}, the Jacobian matrix is explicitly computed and right-preconditioned at every iteration. However, when the RFM is applied to three-dimensional nonlinear PDEs, the resulting residual system $\mathbf{F}(\mathbf{u})$ is typically very large and structurally complex, making explicit Jacobian construction prohibitively expensive. Numerical experiments further show that the dominant computational cost of the INP method is concentrated in the preconditioning stage. To address this issue, we propose the adaptive multi-step inexact preconditioned Newton method, summarized in Algorithm~\ref{alg:amipn}, which requires only one Jacobian evaluation and preconditioning in each outer iteration. The resulting preconditioned Jacobian is reused across multiple approximate Newton steps in the inner loop. In addition, AMIPN prescribes a maximum number of inner iterations together with an adaptive early-stopping criterion to avoid redundant iterations and determine whether the current Jacobian and its preconditioner need to be updated. Moreover, a derivative-free line search is incorporated within the inner loop to ensure monotone residual reduction.

Specifically, at the $k$-th outer iteration, the first inner step is obtained by
\begin{equation}\label{eq:amipn1}
\mathbf{d}_{k,0} = \mathbf{R}_k^{-1} \text{LSQR}\left(\tilde{\mathbf{J}}_k,-\mathbf{F}(\mathbf{u}_k);\eta\right),
\end{equation}
with the update
\begin{align*}
    \mathbf{z}_{k, 1}=\mathbf{u}_k+\alpha_0 \mathbf{d}_{k, 0}, \quad \mathbf{F}_{k, 1}=-\mathbf{F}\left(\mathbf{z}_{k, 1}\right),
\end{align*}
where $\alpha_0>0$ is selected by a derivative-free line search to ensure sufficient decrease of the residual norm. Subsequent approximate Newton steps are computed as
\begin{equation}\label{eq:amipn2}
\mathbf{d}_{k,i} = \mathbf{R}_k^{-1} \text{LSQR}\left(\tilde{\mathbf{J}}_k,-\mathbf{F}(\mathbf{z}_{k,i});\eta\right),
\quad i=1,\ldots,i_k-1,
\end{equation}
together with the updates
\begin{align*}
    \mathbf{z}_{k, i+1}=\mathbf{z}_{k, i}+\alpha_i \mathbf{d}_{k, i}, \quad \mathbf{F}_{k, i+1}=-\mathbf{F}\left(\mathbf{z}_{k, i+1}\right) .
\end{align*}
The actual number of inner steps $i_k$ is adaptively determined by the early-stopping criterion
\begin{align*} 
    |\left\|\mathbf{F}_{k, i+1}-\mathbf{F}_{k, i}\right\||<\tau_{\mathrm{rel}}\left\|\mathbf{F}_{k, 0}\right\|,
 \end{align*}
where $\tau_{\mathrm{rel}}$ is a tolerance. If satisfied, the inner loop terminates early and the current preconditioned Jacobian is retained; otherwise, the process continues until the maximum inner iterations $m_{\max}$ are reached, after which a new preconditioned Jacobian and preconditioner are recomputed. The resulting trial step of the $k$-th outer iteration is
\begin{align*}
\mathbf{s}_k = \sum_{i=0}^{i_k-1}\alpha_i\mathbf{d}_{k,i}.
\end{align*}
These mechanisms reduce the frequency of Jacobian evaluations and preconditioning, lower the overall computational cost, and preserve robustness and accuracy, thus yielding higher efficiency than IPN for large-scale, ill-conditioned nonlinear least-squares problems from nonlinear PDE discretizations.

\begin{algorithm}[!htbp]
\caption{Adaptive multi-step inexact right-preconditioned Newton (AMIPN)}
\label{alg:amipn}
\begin{algorithmic}[1]
\Require $\mathbf{F}(\mathbf{u})$, $\mathbf{J}(\mathbf{u})$, initial guess $\mathbf{x}_0$, tolerances $\varepsilon$ (outer) and $\eta$ (LSQR),
maximum outer iterations $K$, maximum inner steps $m_{\max}$, inner loop early-stopping threshold $\tau_{\text{rel}}$
\Ensure approximate solution $\mathbf{u}$
\State Form $\mathbf{J}_0 \leftarrow \mathbf{J}\left(\mathbf{x}_0\right)$ and perform Steps 2-3 of the IPN method (Algorithm~\ref{alg:irpn}) to obtain the preconditioned Jacobian $\tilde{\mathbf{J}}_0$ and the right preconditioner $\mathbf{R}^{-1}_0$.
\For{$k=0,1,\dots,K-1$}
  \State Set $\mathbf{z}_{k,0}\gets \mathbf{u}_k$, $\mathbf{F}_{k,0}\gets -\mathbf{F}(\mathbf{u}_k)$, $m_k \gets m_{\max}$, \textsc{Flag}$\gets$ True.
  \For{$i=0$ \textbf{to} $m_k-1$}
    \State Compute the inner direction $\mathbf{p}_{k,i}\gets \text{LSQR}(\tilde{\mathbf{J}}_k, \mathbf{F}_{k,i};\,\eta)$.
    \State Compute the inner inexact Newton direction $\mathbf{d}_{k,i}\gets \mathbf{R}_k^{-1} \mathbf{p}_{k,i}$.
    \State choose $\alpha_i>0$ by golden-section search such that $\|\mathbf{F}(\mathbf{z}_{k,i}+\alpha_i \mathbf{d}_{k,i})\|$ sufficiently decreases.
    \State update $\mathbf{z}_{k,i+1}\gets \mathbf{z}_{k,i}+\alpha_i \mathbf{d}_{k,i}$ and $\mathbf{F}_{k,i+1}\gets -\mathbf{F}(\mathbf{z}_{k,i+1})$.
    \If{$\|\mathbf{F}_{k,i+1}-\mathbf{F}_{k,i}\|<\tau_{\text{rel}}\|\mathbf{F}_{k,0}\|$}
      \State \textsc{Flag}$\gets$ False \; \textbf{break}.
    \EndIf
  \EndFor
  \State Update the next iterate as $\mathbf{u}_{k+1}\gets \mathbf{z}_{k,i+1}$ and $\mathbf{F}_{k+1}\gets \mathbf{F}(\mathbf{u}_{k+1})$.
  \If{$|\|\mathbf{F}(\mathbf{u}_{k+1})\|-\|\mathbf{F}(\mathbf{u}_k)\||<\varepsilon$}
    \State \textbf{break}.
  \EndIf
  \If{\textsc{Flag}}
    \State Form $\mathbf{J}_{k+1} \leftarrow \mathbf{J}(\mathbf{u}_{k+1})$ and perform Steps 2-3 of the IPN method (Algorithm~\ref{alg:irpn}) to obtain the preconditioned Jacobian $\tilde{\mathbf{J}}_{k+1}$ and the right preconditioner $\mathbf{R}^{-1}_{k+1}$.
  \EndIf
\EndFor
\State \Return $\mathbf{u} \gets \mathbf{u}_{k+1}$.
\end{algorithmic}
\end{algorithm}

\begin{remark}
In Algorithms~\ref{alg:irpn} and \ref{alg:amipn}, employing the golden-section line search offers a substantial advantage in memory efficiency. This approach requires only multiple evaluations of the objective function $\|\mathbf{F}(\mathbf{x}_k+\alpha_k \mathbf{d}_k)\|$ without storing any additional preconditioned matrix. Once the Jacobian $\mathbf{J}_k$ and its preconditioner $\mathbf{R}$ are computed, updates can be applied in place on $\mathbf{J}_k$ to form the preconditioned Jacobian, thereby reducing the memory footprint by nearly half while maintaining comparable search efficiency. This property is particularly beneficial in large-scale computations.

In contrast, using Armijo, Wolfe, or strong Wolfe line searches requires repeatedly accessing $\mathbf{J}_k^{T}$ and $\mathbf{F}_k$ during each line search at the $k$-th iteration, which not only increases the computational cost but also necessitates storing both the Jacobian and its preconditioned form, thereby leading to substantially higher memory consumption.
\end{remark}

\section{RFM equipped with the proposed nonlinear solvers}\label{sec:5}
In this section, we discretize the nonlinear PDE \eqref{boundary-value problem} using the RFM and combine it with the nonlinear solver introduced in the previous section to develop a complete numerical solution algorithm. The computational domain is defined as
$
\Omega=\{(x,y,z)\mid x\in[a_1,b_1],\; y\in[a_2,b_2],\; z\in[a_3,b_3]\}.
$
The domain $\Omega$ is partitioned—either uniformly or non-uniformly—into $N_x$, $N_y$, and $N_z$ subdomains along the $x$-, $y$-, and $z$-directions, respectively, yielding $M_p=N_xN_yN_z$ subdomains in total. Each subdomain is denoted by
$
\Omega_{e_{mnl}}=[X_m,X_{m+1}]\times[Y_n,Y_{n+1}]\times[Z_l,Z_{l+1}],
$
where $0\le m\le N_x-1$, $0\le n\le N_y-1$, and $0\le l\le N_z-1$, and $e_{mnl}$ is the subdomain index. On each subdomain $\Omega_{e_{mnl}}$, collocation points are placed as
$
(x_p^{e_{mnl}},\, y_q^{e_{mnl}},\, z_r^{e_{mnl}}),\;
0\le p\le Q_x-1,\;0\le q\le Q_y-1,\;0\le r\le Q_z-1,
$
where $x_p^{e_{mnl}}, y_q^{e_{mnl}}, z_r^{e_{mnl}}$ are sets of $Q_x$, $Q_y$, and $Q_z$ points distributed on $[X_m,X_{m+1}]$, $[Y_n,Y_{n+1}]$, and $[Z_l,Z_{l+1}]$, respectively. These points can, for example, be uniformly placed within each subinterval including the endpoints.

We employ the affine mapping introduced in Section~\ref{sec:2.2} to transform the collocation points on each subdomain $\Omega_{e_{mnl}}$ into the reference cube $[-1,1]^3$. Using the PoU functions $\psi_i^a(\mathbf{x})$ defined in \eqref{RFM_approximate_sol}, the global RFM solution is represented in a piecewise manner over the $N_xN_yN_z$ subdomains. On each subdomain $\Omega_{e_{mnl}}$, the local solution is expressed as a linear combination of $J$ random feature functions, with the same number of basis functions used on all subdomains for consistency. The resulting global approximation takes the form
\begin{align*}
    u(x,y,z) =
\begin{cases}
\displaystyle \sum_{j=1}^J \phi_j^{000}(x,y,z)\,u_j^{000}, & (x,y,z)\in \Omega_{e_{000}}, \\[2ex]
\displaystyle \sum_{j=1}^J \phi_j^{001}(x,y,z)\,u_j^{001}, & (x,y,z)\in \Omega_{e_{001}}, \\[1ex]
\qquad \vdots & \qquad \vdots \\[1ex]
\displaystyle \sum_{j=1}^J \phi_j^{mnl}(x,y,z)\,u_j^{mnl}, & (x,y,z)\in \Omega_{e_{mnl}}, \\[1ex]
\qquad \vdots & \qquad \vdots \\[1ex]
\displaystyle \sum_{j=1}^J \phi_j^{N_x-1,\,N_y-1,\,N_z-1}(x,y,z)\,u_j^{N_x-1,\,N_y-1,\,N_z-1}, & (x,y,z)\in \Omega_{e_{N_x-1,\,N_y-1,\,N_z-1}} .
\end{cases}
\end{align*}
On each subdomain $\Omega_{e_{mnl}}$, the local approximation takes the form
\begin{equation}\label{local_sol}
    u^{mnl}(x,y,z)=\sum_{j=1}^J \phi_j^{mnl}(x,y,z)\,u_j^{mnl}.
\end{equation}
Substituting \eqref{local_sol} into the governing equation and enforcing it at the collocation point 
$(x_p^{e_{mnl}},y_q^{e_{mnl}},z_r^{e_{mnl}})$ yields
\begin{equation}\label{inner_nonlinear_system}
\lambda\sum_{j=1}^J\big(\mathcal{L}\phi_j^{mnl}\big)\Big|_{(x_p^{e_{mnl}},y_q^{e_{mnl}},z_r^{e_{mnl}})} u_j^{mnl}
+\mu\,\mathcal{N}\!\Big(\sum_{j=1}^J\phi_j^{mnl}u_j^{mnl}\Big)\Big|_{(x_p^{e_{mnl}},y_q^{e_{mnl}},z_r^{e_{mnl}})}
=f(x_p^{e_{mnl}},y_q^{e_{mnl}},z_r^{e_{mnl}}).
\end{equation}
where $0\le m\le N_x-1,\;0\le n\le N_y-1,\;0\le l\le N_z-1,\;0\le p\le Q_x-1,\;0\le q\le Q_y-1,\;0\le r\le Q_z-1$.

We assume that in the nonlinear PDE \eqref{boundary-value problem}, at least one of the 
operators $\mathcal{L}$ or $\mathcal{N}$ is of order $k$. 
On the interfaces $x=X_{m+1}$ ($0\le m\le N_x-2$), the $C^0,\dots,C^{k-1}$ continuity reads
\begin{equation}\label{con:1}
\left\{
\begin{aligned}
&\sum_{j=1}^J \frac{\partial^{\,s} \phi_j^{mnl}}{\partial x^s}\!\left(X_{m+1}, y_q^{mnl}, z_r^{mnl}\right) u_j^{mnl}
-\sum_{j=1}^J \frac{\partial^{\,s} \phi_j^{m+1,n,l}}{\partial x^s}\!\left(X_{m+1}, y_q^{m+1,n,l}, z_r^{m+1,n,l}\right) u_j^{m+1,n,l} = 0, 
\quad s=0,1,\dots,k-1, \\
&\text{for } 0 \leq m \leq N_x-2,\; 0 \leq n \leq N_y-1,\; 0 \leq l \leq N_z-1,\; 0 \leq q \leq Q_y-1,\; 0 \leq r \leq Q_z-1 .
\end{aligned}
\right.
\end{equation}
Similarly, on $y=Y_{n+1}$ ($0\le n\le N_y-2$),
\begin{equation}\label{con:2}
\left\{
\begin{aligned}
&\sum_{j=1}^J \frac{\partial^{\,s} \phi_j^{mnl}}{\partial y^s}\!\left(x_p^{mnl}, Y_{n+1}, z_r^{mnl}\right) u_j^{mnl}
-\sum_{j=1}^J \frac{\partial^{\,s} \phi_j^{m,n+1,l}}{\partial y^s}\!\left(x_p^{m,n+1,l}, Y_{n+1}, z_r^{m,n+1,l}\right) u_j^{m,n+1,l} = 0, 
\quad s=0,1,\dots,k-1,\\
&\text{for } 0 \leq m \leq N_x-1,\; 0 \leq n \leq N_y-2,\; 0 \leq l \leq N_z-1,\; 0 \leq p \leq Q_x-1,\; 0 \leq r \leq Q_z-1 .
\end{aligned}
\right.
\end{equation}
On $z=Z_{l+1}$ ($0\le l\le N_z-2$),
\begin{equation}\label{con:3}
\left\{
\begin{aligned}
&\sum_{j=1}^J \frac{\partial^{\,s} \phi_j^{mnl}}{\partial z^s}\!\left(x_p^{mnl}, y_q^{mnl}, Z_{l+1}\right) u_j^{mnl}
-\sum_{j=1}^J \frac{\partial^{\,s} \phi_j^{m,n,l+1}}{\partial z^s}\!\left(x_p^{m,n,l+1}, y_q^{m,n,l+1}, Z_{l+1}\right) u_j^{m,n,l+1} = 0, 
\quad s=0,1,\dots,k-1,\\
&\text{for } 0 \leq m \leq N_x-1,\; 0 \leq n \leq N_y-1,\; 0 \leq l \leq N_z-2,\; 0 \leq p \leq Q_x-1,\; 0 \leq q \leq Q_y-1 .
\end{aligned}
\right.
\end{equation}

Applying the boundary condition \eqref{boundary-value problem} to the spatial boundaries 
$x=a_1, b_1$, $y=a_2, b_2$, and $z=a_3, b_3$ yields
\begin{equation}\label{bdy:1}
\left\{
\begin{aligned}
&\sum_{j=1}^J \mathcal{B}\!\left(\phi_j^{mnl}(a_1, y_q^{mnl}, z_r^{mnl})\right) u_j^{mnl}
- g(a_1, y_q^{mnl}, z_r^{mnl}) = 0, \\[1ex]
&\sum_{j=1}^J \mathcal{B}\!\left(\phi_j^{mnl}(b_1, y_q^{mnl}, z_r^{mnl})\right) u_j^{mnl}
- g(b_1, y_q^{mnl}, z_r^{mnl}) = 0, \\[1ex]
&\text{for } 0 \leq n \leq N_y-1,\; 0 \leq l \leq N_z-1,\; 0 \leq q \leq Q_y-1,\; 
0 \leq r \leq Q_z-1,\; m=0\,(x=a_1) \text{ or } m=N_x-1\,(x=b_1);
\end{aligned}
\right.
\end{equation}
\begin{equation}\label{bdy:2}
\left\{
\begin{aligned}
&\sum_{j=1}^J \mathcal{B}\!\left(\phi_j^{mnl}(x_p^{mnl}, a_2, z_r^{mnl})\right) u_j^{mnl}
- g(x_p^{mnl}, a_2, z_r^{mnl}) = 0, \\[1ex]
&\sum_{j=1}^J \mathcal{B}\!\left(\phi_j^{mnl}(x_p^{mnl}, b_2, z_r^{mnl})\right) u_j^{mnl}
- g(x_p^{mnl}, b_2, z_r^{mnl}) = 0,  \\[1ex]
&\text{for } 0 \leq m \leq N_x-1,\; 0 \leq l \leq N_z-1,\; 0 \leq p \leq Q_x-1,\; 0 \leq r \leq Q_z-1,
\; n=0(y=a_2) \text{ or } n=N_y-1(y=b_2);
\end{aligned}
\right.
\end{equation}
\begin{equation}\label{bdy:3}
\left\{
\begin{aligned}
&\sum_{j=1}^J \mathcal{B}\!\left(\phi_j^{mn0}(x_p^{mn0}, y_q^{mn0}, a_3)\right) u_j^{mn0}
- g(x_p^{mn0}, y_q^{mn0}, a_3) = 0, \\[1ex]
&\sum_{j=1}^J \mathcal{B}\!\left(\phi_j^{mnl}(x_p^{mnl}, y_q^{mnl}, b_3)\right) u_j^{mnl}
- g(x_p^{mnl}, y_q^{mnl}, b_3) = 0,  \\[1ex]
&\text{for } 0 \leq m \leq N_x-1,\; 0 \leq n \leq N_y-1,\; 0 \leq p \leq Q_x-1,\; 0 \leq q \leq Q_y-1,
\; l=0(z=a_3) \text{ or } l=N_z-1(z=b_3).
\end{aligned}
\right.
\end{equation}

Building on the discretization of \eqref{boundary-value problem}, the overall RFM framework—combined with the nonlinear solver, either the INP or AMIPN method—is summarized in the following algorithm.

\begin{algorithm}
\caption{RFM with IPN/AMIPN nonlinear solvers}
\label{alg:RFM3D}
\begin{algorithmic}[1]
\Require Partition $\{\Omega_{e_{mnl}}\}$; collocation sets 
$\{x_p^{e_{mnl}}\},\{y_q^{e_{mnl}}\},\{z_r^{e_{mnl}}\}$; 
basis functions $\{\phi_j^{mnl}\}_{j=1}^J$; 
operators $\mathcal L,\mathcal N,\mathcal B$; 
data $f,g$; continuity order $k$; 
initial guess $\mathbf{u}_0$; tolerance $\varepsilon$.
\Ensure Approximate solution.

\State \textbf{Step 1: Interior assembly.} The nonlinear residual of the governing equations, derived from \eqref{inner_nonlinear_system}, is given by 
\begin{align*}
    \mathbf{R}_{pde}(\mathbf{u})=\lambda \mathbf{A}_{\mathrm{int}} \mathbf{u}+\mu \mathbf{g}_{\mathrm{int}}(\mathbf{u})-\mathbf{f}_{\mathrm{int}}.
\end{align*}
\State \textbf{Step 2: Continuity assembly.} The continuity residual, assembled from \eqref{con:1}, \eqref{con:2}, and \eqref{con:3}, is
\begin{align*}
    \mathbf{R}_{\mathrm{con}}(\mathbf{u})=A_{\mathrm{con}} \mathbf{u}.
\end{align*}
\State \textbf{Step 3: Boundary conditions.} The boundary residual, assembled from \eqref{bdy:1}, \eqref{bdy:2}, and \eqref{bdy:3}, is given by
\begin{align*}
    \mathbf{R}_{\mathrm{bdy}}(\mathbf{u})=\mathbf{A}_{\mathrm{bdy}} \mathbf{u}-\mathbf{b}_{\mathrm{bdy}} .
\end{align*}
\State \textbf{Step 4: Global nonlinear system.} Define
\begin{equation}\label{equation:5.13}
\mathbf{F}(\mathbf{u})=
\begin{bmatrix}
\mathbf{R}_{pde}(\mathbf{u})\\
\mathbf{R}_{\mathrm{con}}(\mathbf{u})\\
\mathbf{R}_{\mathrm{bdy}}(\mathbf{u})
\end{bmatrix},
\quad
\min_{\mathbf{u}} \tfrac{1}{2}\|\mathbf{F}(\mathbf{u})\|_2^2
\quad \text{equivalently, solve } \mathbf{F}(\mathbf{u})=0.
\end{equation}

\State \textbf{Step 5: Nonlinear solve.} Call the nonlinear solver (Algorithm \ref{alg:irpn} or Algorithm \ref{alg:amipn}) to solve \eqref{equation:5.13}.
\State \Return $\mathbf{u}_k$, the approximate solution.
\end{algorithmic}
\end{algorithm}

Here, $\mathbf{A}_{\text{int}}$ is obtained by applying the linear operator $\mathcal{L}$ to the random feature basis functions, $\mathbf{g}_{\text{int}}(\mathbf{u})$ represents the discrete nonlinear contribution, and $\mathbf{f}_{\text{int}}$ denotes the source term vector. The matrices $\mathbf{A}_{\text{con}}$ and $\mathbf{A}_{\text{bdy}}$ correspond to the continuity and boundary constraints, respectively, with $\mathbf{b}{\text{bdy}}$ denoting the discrete boundary data. The Jacobian of $\mathbf{F}(\mathbf{u})$ is then
\[
J(\mathbf{u}) = \frac{\partial \mathbf{F}(\mathbf{u})}{\partial \mathbf{u}}
=\begin{bmatrix}
\lambda\mathbf{A}_{\mathrm{int}}+\mu\frac{\partial \mathbf{g}_{\mathrm{int}}(\mathbf{u})}{\partial \mathbf{u}}\\[3pt]
\mathbf{A}_{\mathrm{con}} \\[3pt]
\mathbf{A}_{\mathrm{bdy}}
\end{bmatrix}.
\]
At the $k$-th Newton iteration, the correction $\delta \mathbf{u}$ is obtained from the least-squares problem
\begin{equation}\label{newton_LS_system}
\min_{\delta\mathbf{u}\in\mathbb{R}^n}
\|\mathbf{F}(\mathbf{u}_k)+\mathbf{J}(\mathbf{u}_k)\delta\mathbf{u}\|_2.
\end{equation}

Since the PDE, continuity, and boundary contributions may differ significantly in magnitude, the system \eqref{newton_LS_system} is rescaled before applying IPN or AMIPN. Specifically, the Jacobian
$
\mathbf{J}_k\left(\mathbf{u}_k\right)=\left(u_{i j}^k\right) \in \mathbb{R}^{m \times n},
$
and residual
$
\mathbf{F}\left(\mathbf{u}_k\right)=\left[f_1\left(\mathbf{u}_k\right), \ldots, f_m\left(\mathbf{u}_k\right)\right]^T \in \mathbb{R}^m,
$
are scaled row-wise. This procedure follows the idea of adaptive weighting: scaling factors are adjusted according to the magnitude of each residual to prevent any single equation from dominating the iteration, thereby improving numerical comparability and stability.

The scaling factor for the $i$-th row is defined as
\begin{align*}
    \lambda_{k, i}=\frac{c}{\max _{1 \leq j \leq n}|u_{i j}^k|}, \quad i=1, \ldots, m,
\end{align*}
with $c=100$ in all experiments. The scaled system is then
\begin{align*}
    u_{i j}^k \leftarrow \lambda_{k, i} u_{i j}^k, \quad f_i \leftarrow \lambda_{k, i} f_i, \quad i=1, \ldots, m, j=1, \ldots, n .
\end{align*}
After this normalization, the system is successively preconditioned and solved.

\section{Numerical experiments}\label{sec:6}
This section presents numerical experiments on two categories of nonlinear PDEs. 
The first set consists of three-dimensional steady-state problems used to assess the numerical convergence of the RFM combined with the proposed nonlinear solvers, while the second set includes two-dimensional time-dependent problems with complex dynamics or intricate geometries to demonstrate their effectiveness and practicality. 
The computational environment and notation follow Section~\ref{subsec:3.4}, with ``NJ'' denoting the number of Jacobian evaluations. 
In Algorithms~\ref{alg:irpn} (IPN) and \ref{alg:amipn} (AMIPN), the sampling factor is fixed at $\gamma=3$, the LSQR termination tolerance is set to $\eta=10^{-6}$, and the overall termination tolerance is prescribed as $\varepsilon=10^{-10}$. 
For AMIPN, the inner-loop termination criterion is chosen as $\tau_{\mathrm{rel}}=10^{-3}$, and the maximum number of inner steps is set to $m_{\max}=3$.


\subsection{Comparison with FEM, FDM, PINN, and WAN}
In this section, we evaluate the accuracy of the proposed nonlinear solvers by comparing them with FDM, FEM, PINN, and WAN on the following three-dimensional nonlinear elliptic problem defined on the unit cube domain $(0,1)^3$ :
\begin{equation}\label{3D_nonlinear_PDE} 
    \begin{cases}-\Delta u + u^3 =f(x, y, z), & (x, y, z) \in (0,1)^3 ,
    \\ u(x, y, z)=g(x, y, z) ,
        & \end{cases} 
\end{equation} 
where the exact solution is prescribed as
\begin{align*}
    u(x, y, z)=\sin (\pi x) \sin (\pi y) \sin (\pi z) .
\end{align*}
Dirichlet boundary conditions are imposed on $\partial(0,1)^3$. 

\begin{figure}[!htbp]
    \renewcommand\figurename{Figure}
    \centering
    \begin{subfigure}[b]{0.40\textwidth}
        \centering
        \includegraphics[width=\textwidth]{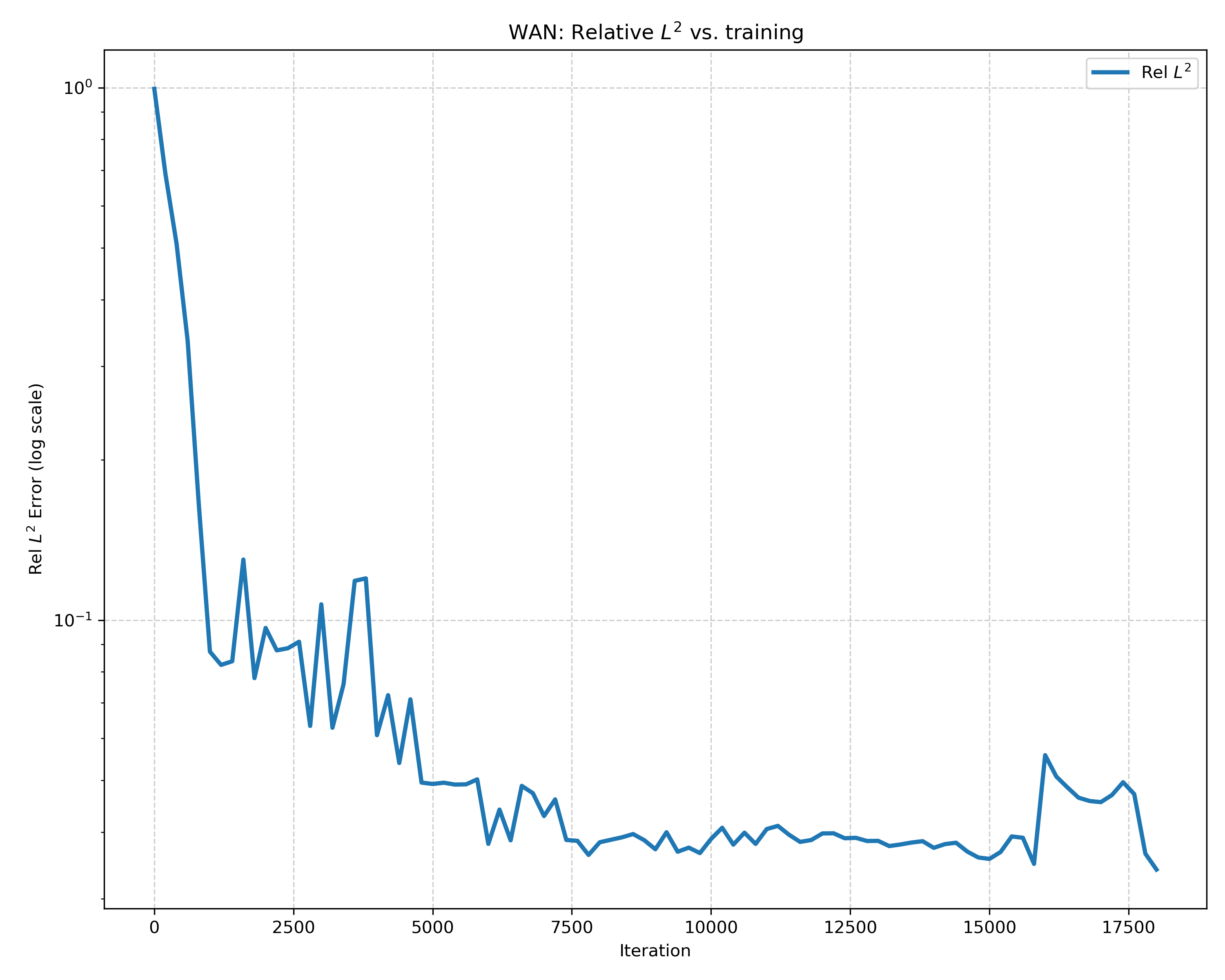}
        \label{fig:wan_l2}
    \end{subfigure}
    \begin{subfigure}[b]{0.40\textwidth}
        \centering
        \includegraphics[width=\textwidth]{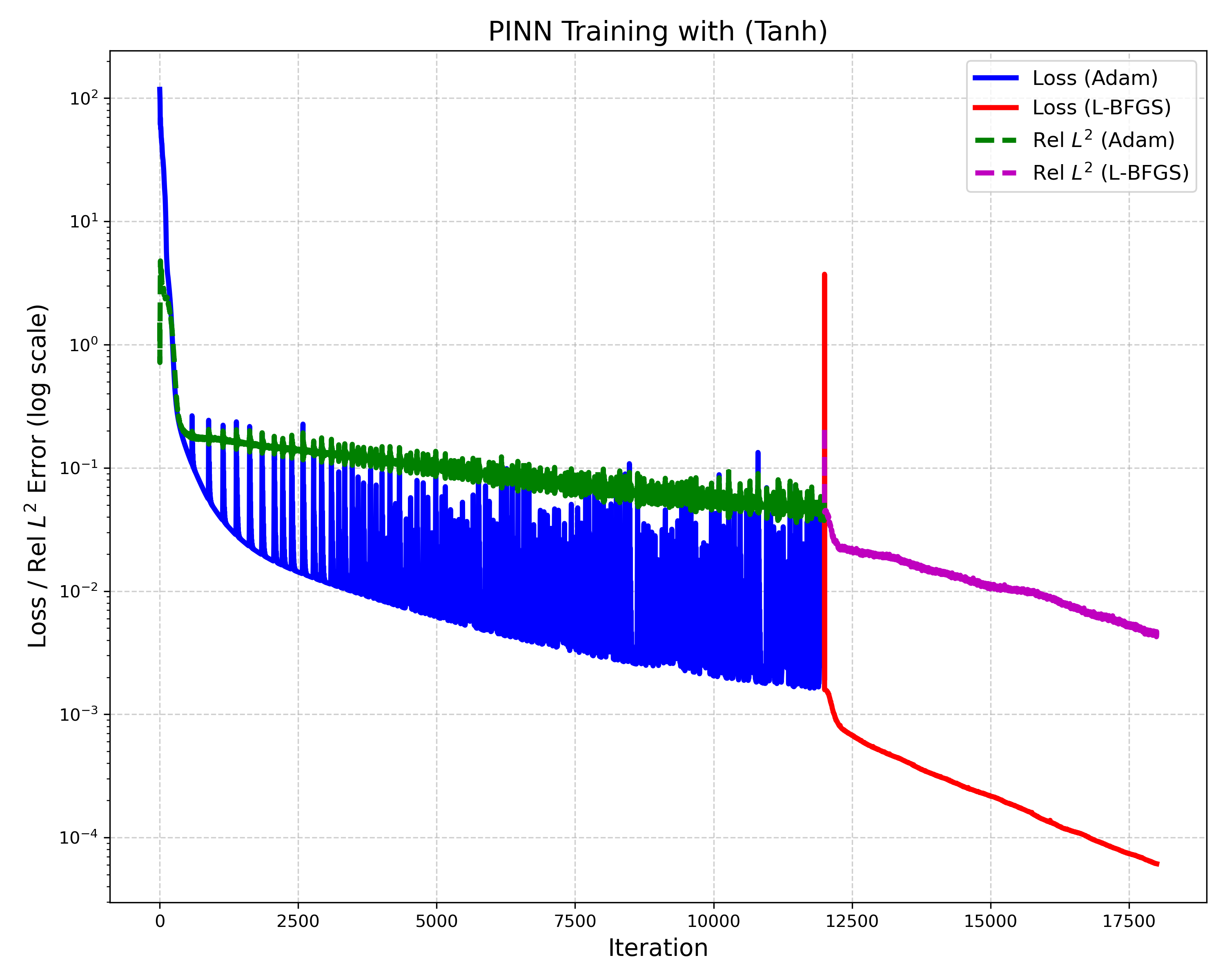}
        \label{fig:3D_nonlinear_cubic}
    \end{subfigure}
    \caption{Left: WAN (Relative $L^2$ vs. training); 
             Right: PINN (Relative $L^2$ and loss vs. training).}
    \label{fig:comparison_2plots}
\end{figure}

Table~\ref{3D_nonlinear_result_one} reports the numerical results of AMIPN, FDM, and FEM, while Figure~\ref{fig:comparison_2plots} presents those of PINN and WAN. Both FEM and FDM are performed on uniform grids of sizes $40^3, 50^3$, and $60^3$. The FEM is implemented using the open-source finite element package FEniCS \cite{alnaes2015fenics} with standard $P_1$ Lagrange elements on uniform tetrahedral meshes. For the FDM, a standard second-order central-difference scheme on uniform Cartesian grids is adopted. The resulting nonlinear algebraic systems are solved by Picard iteration with a stopping criterion that the relative change between successive iterates is below $10^{-6}$. The PINN employs a fully connected feedforward architecture with five hidden layers of 40 neurons each, whereas the WAN consists of a solution network and an adversarial test network, each with five hidden layers of 50 neurons. The training datasets comprise 18,000 interior points and 6,000 boundary points, both generated via Latin hypercube sampling.

Table~\ref{3D_nonlinear_result_one} shows that, with the same number of collocation points ($\sqrt[3]{Q}=20,25,30$), the AMIPN method produces errors 3–7 orders of magnitude smaller than those of FDM and FEM. The superiority of AMIPN becomes even more evident as both collocation points and basis functions increase. For instance, at $\sqrt[3]{Q}=25$, the errors of FEM and FDM are approximately $1.79\times 10^{-3}$ and $3.29\times 10^{-4}$, while AMIPN achieves only $3.12\times 10^{-9}$ with $J=1600$. Moreover, for a fixed number of collocation points, enlarging the basis set $J$ yields further accuracy gains: when $\sqrt[3]{Q}=30$, the error decreases from $3.05\times 10^{-9}$ ($J=1600$) to $6.20\times 10^{-11}$ ($J=2400$), demonstrating excellent convergence. Figure~\ref{fig:comparison_2plots} confirms this trend, showing that for problem~\eqref{3D_nonlinear_PDE} the final relative $L^2$ errors of PINN and WAN remain at $O(10^{-3})$ and $O(10^{-2})$, respectively—several orders of magnitude larger than those of AMIPN.

\begin{table}[!htbp]
  \centering
  \caption{Error comparison of AMIPN, FDM and FEM for solving the 3D nonlinear elliptic equation \eqref{3D_nonlinear_PDE}.}
  \label{3D_nonlinear_result_one}
  \begin{tabular}{c c c c c c c}
    \toprule
    $(N_x, N_y, N_z)$ & $\sqrt[3]{Q} $ & $J$ & Solver & IT & NJ & $u$ error \\
    \midrule
    \multirow{12}{*}{$(2, 2, 2)$} 
    & \multirow{4}{*}{$20$}  & \multirow{1}{*}{400} & AMIPN & $2$ & $2$ & $ 5.38e-05$ \\
        && \multirow{1}{*}{800} & AMIPN & $3$ & $2$& $ 8.58e-07$ \\
        \cmidrule(lr){3-7}
        &&& FEM & $4$ & \multicolumn{1}{c}{$--$}& $2.80e-03$ \\
        &&& FDM & $6$ & \multicolumn{1}{c}{$--$}& $5.19e-04$ \\
    \cmidrule(lr){2-7}
        & \multirow{4}{*}{$25$} 
        &\multirow{1}{*}{1200}& AMIPN & $3$ & $2$ & $3.96e-08$ \\
        &&\multirow{1}{*}{1600}& AMIPN & $3$ & $2$& $ 3.12e-09$ \\
        \cmidrule(lr){3-7}
        &&& FEM & $4$ & \multicolumn{1}{c}{$--$}& $1.79e-03$ \\
        &&& FDM & $6$ & \multicolumn{1}{c}{$--$}& $3.29e-04 $ \\
    \cmidrule(lr){2-7}
        & \multirow{5}{*}{$30$} 
        &\multirow{1}{*}{1600}& AMIPN & $3$ & $2$ & $3.05e-09$ \\
         &&\multirow{1}{*}{2000}& AMIPN & $3$ & $2$ & $3.66e-10$ \\
         &&\multirow{1}{*}{2400}& AMIPN & $3$ & $2$ & $6.20e-11$ \\ 
         \cmidrule(lr){3-7}
        &&& FEM & $4$ & \multicolumn{1}{c}{$--$}& $1.25e-03$ \\
        &&& FDM & $6$ & \multicolumn{1}{c}{$--$}& $2.27e-04$ \\
    \bottomrule
  \end{tabular}
\end{table}

\subsection{Three-dimensional steady nonlinear PDEs}
In this section, we consider a set of three-dimensional steady-state partial differential equations—including a strongly nonlinear elliptic equation, a nonlinear diffusion-reaction equation, a nonlinear Helmholtz equation, and the Gray-Scott system—to offer a more comprehensive evaluation of the proposed methods.

To enable a direct comparison with the numerical results in Section~\ref{subsec:3.4}, 
we first apply the AMIPN method to the three-dimensional nonlinear elliptic problem 
\eqref{benchmark_nonlinear_elliptic_equation}; the results are summarized in 
Table~\ref{3D_nonlinear_equation}. The table reports the time required to assemble 
the NLS system via automatic differentiation, the number of outer iterations (IT) 
and Jacobian evaluations (NJ), the CPU time for solving the NLS system, the final 
residual norm $\|\mathbf{F}(\mathbf{u}_k)\|_2$, and the relative $L^2$ and $H^1$ errors. 

With $(N_x,N_y,N_z)=(2,2,2)$, increasing $\sqrt[3]{Q}$ and $J$ steadily improves accuracy: 
the $L^2$ error decreases from $O(10^{-6})$ to $O(10^{-11})$, 
the $H^1$ error from $O(10^{-5})$ to $O(10^{-9})$, 
and the residual norm from $O(10^{-4})$ to $O(10^{-11})$. 
Meanwhile, the number of outer iterations stabilizes at five, with four Jacobian 
evaluations per iteration, whereas the assembly time and CPU time for solving the 
NLS system increase monotonically with problem size—from about 23\,s to 410\,s and 
from 201\,s to 2745\,s, respectively. These results, together with those in 
Tables~\ref{3D_nonlinear_equation_newton_LM} and \ref{3D_nonlinear_equation}, 
demonstrate that the AMIPN method achieves substantially higher accuracy and efficiency 
than classical solvers such as Newton, LM, Mor{\'e}, and TR, particularly for large-scale 
nonlinear elliptic problems.

To further assess the accuracy of the proposed method, 
we consider the following strongly nonlinear elliptic equation posed on the unit cube $\Omega=(0,1)^3$:
\begin{equation}\label{strongly_nonlinear_elliptic_PDE}
\begin{cases}
-\nabla \cdot(\Lambda(\nabla u) \nabla u)=f(x,y,z), & (x,y,z)\in \Omega, \\[3pt]
u(x,y,z)=g(x,y,z), & (x,y,z)\in \partial\Omega,
\end{cases}
\end{equation}
where $\Lambda(\nabla u)=1+|\nabla u|^2$. 
The exact solution is
\[
u(x,y,z)=\sin(x(1-x))\sin(y(1-y))e^{z+1/2},
\]
with Dirichlet boundary conditions prescribed on $\partial\Omega$.

Table~\ref{3D_PDE_results} reports the detailed numerical results for Example~\eqref{strongly_nonlinear_elliptic_PDE}, 
including the relative $L^2$ errors of $u$ and its first-order derivatives, 
the final residual norm $\|\mathbf{F}(\mathbf{u}_k)\|_2$, the number of outer iterations, 
and the actual number of Jacobian evaluations. 
As $J$ increases from 400 to 2400, the $L^2$ error of $u$ decreases from $O(10^{-5})$ to $O(10^{-11})$, 
and the errors of $u_x$, $u_y$, and $u_z$ decrease from $O(10^{-3}\!-\!10^{-4})$ to $O(10^{-9})$, 
exhibiting an approximately exponential convergence trend. 
These results demonstrate that the AMIPN method achieves extremely high accuracy 
while converging in only a few iterations, even for this three-dimensional strongly nonlinear elliptic PDE.

\begin{table}[!htbp]
  \centering
  \caption{The results of  the AMIPN method for the 3D nonlinear elliptic problem \eqref{benchmark_nonlinear_elliptic_equation}}
  \label{3D_nonlinear_equation}
  \begin{tabular}{c c c c c c c c c c c}
    \toprule
    $(N_x, N_y, N_z)$ & $\sqrt[3]{Q} $ & $J$ & Solver & IT & NJ & Assemble$(s)$&CPU&$\|\mathbf{F}(\mathbf{u}_k)\|_2$ & $u$ error & $H^1$ error \\
    \midrule
    \multirow{6}{*}{$(2, 2, 2)$} & \multirow{1}{*}{$20$} 
        & \multirow{1}{*}{800} & AMIPN & $6$ & $4$ &$23.76$&$ 201.00$ &$ 1.07e-04$& $1.02e-06$ & $2.85e-05$\\
    \cmidrule(lr){2-11}
        & \multirow{2}{*}{$25$} 
        & \multirow{1}{*}{1200} & AMIPN  & $5$ & $4$ &$69.35$ & $ 580.62$&$3.04e-06$ & $ 3.87e-08$ & $1.32e-06$\\
        && \multirow{1}{*}{1600} & AMIPN  & $5$ & $4$ &$117.31$ & $ 967.20$ &$2.42e-08$ & $ 3.89e-09$ & $ 1.48e-07$\\
      \cmidrule(lr){2-11}
        & \multirow{3}{*}{$30$} 
        & \multirow{1}{*}{1600} & AMIPN  & $5$ & $4$ &$187.94$ & $ 1335.61$&$ 3.96e-08$ & $3.41e-09$ & $1.30e-07$\\
        && \multirow{1}{*}{2000} & AMIPN  & $5$ & $4$ &$ 282.08$ & $2176.57$ &$ 5.96e-10$ & $4.16e-10$ & $1.77e-08$\\
        && \multirow{1}{*}{2400} & AMIPN  & $5$ & $4$ &$ 410.16$ & $ 2744.73$ &$ 1.21e-11$ & $7.04e-11$ &  $3.05e-09$\\
    \bottomrule
  \end{tabular}
\end{table}

\begin{table}[!htbp]
    \centering
    \caption{The results of  the AMIPN method for the 3D strongly nonlinear elliptic problem \eqref{strongly_nonlinear_elliptic_PDE}}
    \label{3D_PDE_results}
    \begin{tabular}{c c c c c c c c c c c}
      \toprule
      $(N_x, N_y, N_z)$ & $\sqrt[3]{Q}$ & $J$ & Solver & IT & NJ &$\|\mathbf{F}(\mathbf{u}_k)\|_2$& $u$ error & $u_x$ error & $u_y$ error & $u_z$ error \\
      \midrule
      \multirow{7}{*}{$(2, 2, 2)$} 
        & \multirow{2}{*}{$20$} & $400$  & AMIPN & $4$ & $3$ &$2.79e-02$& $ 7.28e-05$ & $4.65e-04$ & $4.57e-04$ & $1.61e-03$ \\
        && $800$  & AMIPN & $4$ & $3$ &$ 2.11e-05 $& $ 1.64e-06$ & $1.38e-05$ & $ 1.67e-05$ & $5.65e-05$ \\
      \cmidrule(lr){2-11}
        & \multirow{2}{*}{$25$} 
        & $1200$ & AMIPN & $4$ & $3$  &$4.09e-08$& $5.25e-08$ & $6.02e-07$ & $6.00e-07$ & $1.90e-06$ \\
        && $1600$ & AMIPN & $4$ & $3$  &$4.83e-10$& $4.91e-09$ & $6.11e-08$ & $ 6.35e-08$ & $ 2.28e-07$ \\
      \cmidrule(lr){2-11}
        & \multirow{3}{*}{$30$} 
        & $1600$ & AMIPN & $4$ & $3$  &$ 4.57e-10$& $3.73e-09$ & $4.84e-08$ & $4.75e-08$ & $1.65e-07$ \\
        && $2000$ & AMIPN & $4$ & $3$  &$1.64e-11$& $6.47e-10$ & $ 9.13e-09$ & $9.08e-09$ & $3.21e-08$ \\
        && $2400$ & AMIPN & $3$ & $3$  &$ 3.35e-13 $& $9.33e-11$ & $1.44e-09$ & $1.45e-09$ & $5.14e-09$ \\
      \bottomrule
    \end{tabular}
  \end{table}

To compare the numerical performance of the IPN and AMIPN methods, 
we consider the following three-dimensional nonlinear Helmholtz equation posed on the unit cube $\Omega=(0,1)^3$:
\begin{equation}\label{3D_nonlinear_Helmholtz_equation}
\begin{cases}
u_{xx}+u_{yy}+u_{zz}-100u+10\cosh(u)=f(x,y,z), & (x,y,z)\in\Omega, \\[3pt]
u(x,y,z)=g(x,y,z), & (x,y,z)\in\partial\Omega,
\end{cases}
\end{equation}
with the exact solution
\[
u(x,y,z)=(2+x^2+y^2+z^2)\sin(\pi x)\sin(\pi y)\sin(\pi z),
\]
and corresponding Dirichlet boundary conditions imposed on $\partial\Omega$.

Table~\ref{3D_nonlinear_helmholtz_equation} reports the performance of the two solvers for the 3D nonlinear Helmholtz equation. As the numbers of collocation points and random basis functions increase, both methods exhibit markedly improved accuracy: the relative $L^2$ error decreases from about $10^{-5}$ at $J=400$ to $10^{-11}$ at $J=2400$, while the $H^1$ error decreases from $O(10^{-4})$ to $O(10^{-9})$. Notably, AMIPN attains essentially the same accuracy as IPN while requiring fewer outer iterations (3-4 vs.\ 5-6) and significantly fewer Jacobian evaluations (2 vs.\ 5-6). Consequently, AMIPN achieves about a $1.7\times$ speedup at smaller problem sizes, which further increases to $2.6\times$ at the largest scale. This trend clearly indicates that the relative advantage of AMIPN over IPN becomes more pronounced as the problem size grows.

\begin{table}[!htbp]
  \centering
  \caption{The results of the IPN and  AMIPN methods for the 3D nonlinear Helmholtz equation \eqref{3D_nonlinear_Helmholtz_equation}}
  \label{3D_nonlinear_helmholtz_equation}
  \begin{tabular}{c c c c c c c c c c}
    \toprule
    $(N_x, N_y, N_z)$ & $\sqrt[3]{Q} $ & $J$ & Solver & IT& NJ & CPU & Speedup & $u$ error & $H^1$ error \\
    \midrule
    \multirow{14}{*}{$(2, 2, 2)$} 
        & \multirow{4}{*}{$20$} & \multirow{2}{*}{400} &IPN & $5$ & $5$ & $ 52.30$& \multirow{2}{*}{1.73} & $2.55e-05$ & $ 5.85e-04$ \\
        &&& AMIPN & $3$ & $2$ & $ 30.21$& & $2.54e-05$ & $5.83e-04$ \\
         \cmidrule(lr){3-10}
        && \multirow{2}{*}{800} &IPN & $5$ & $5$ & $170.53$ & \multirow{2}{*}{1.94} & $ 4.33e-07$ & $ 1.31e-05$ \\
        &&& AMIPN & $3$ & $2$ & $87.87$ & & $4.32e-07$ & $1.31e-05$ \\
    \cmidrule(lr){2-10}
        & \multirow{4}{*}{$25$} &\multirow{2}{*}{1200} &IPN & $5$ & $5$ & $ 485.91$ &  \multirow{2}{*}{2.07} & $1.91e-08$ & $7.15e-07$ \\
        &&& AMIPN & $4$ & $2$ & $234.66$& & $1.91e-08$ & $7.14e-07$ \\
        \cmidrule(lr){3-10}
        && \multirow{2}{*}{1600} &IPN & $5$ & $5$ & $ 1074.96$ &  \multirow{2}{*}{2.40} & $1.57e-09$ & $6.78e-08$ \\
        &&& AMIPN & $4$ & $2$ & $437.56$& & $1.57e-09$ & $ 6.78e-08$ \\
      \cmidrule(lr){2-10}
        & \multirow{6}{*}{$30$} & \multirow{2}{*}{1600} &IPN & $5$ & $5$ & $ 1396.71$ &  \multirow{2}{*}{2.35} & $1.48e-09$ & $ 6.17e-08$ \\
        &&& AMIPN & $4$ & $2$ &$593.49$& & $ 1.47e-09$ & $ 6.16e-08$ \\
        \cmidrule(lr){3-10}
        && \multirow{2}{*}{2000} &IPN & $6$ & $6$ & $ 2147.64$& \multirow{2}{*}{2.41} & $ 1.91e-10$ & $8.87e-09$ \\
        &&& AMIPN & $4$ & $2$ &$891.53$ & & $9.48e-09$ & $3.17e-07$ \\
        \cmidrule(lr){3-10}
        && \multirow{2}{*}{2400} &IPN & $6$ & $6$ & $ 3386.93$& \multirow{2}{*}{2.59} & $ 2.64e-11$ & $1.26e-09$ \\
        &&& AMIPN & $4$ & $2$ &$1313.53$ & & $2.56e-11$ & $1.26e-09$ \\
    \bottomrule
  \end{tabular}
\end{table}

We consider the following three-dimensional steady-state nonlinear diffusion–reaction equation posed on the unit ball $\Omega=\{(x,y,z)\in\mathbb{R}^3: x^2+y^2+z^2\leq 1\}$:
\begin{equation}\label{Allen-Cahn-Bratu_equation}
\begin{cases}
-\nabla\cdot(\Lambda(u)\nabla u) + (u^3-u) + e^u = f(x,y,z), & (x,y,z)\in\Omega, \\[3pt]
u(x,y,z)=g(x,y,z), & (x,y,z)\in\partial\Omega,
\end{cases}
\end{equation}
where $\Lambda(u)=1+u^2$. The exact solution is chosen as
\[
u(x,y,z)=(x^2+y^2+z^2)\sin(\pi x)\sin(\pi y)\sin(\pi z),
\]
and $f(x,y,z)$ is determined accordingly so that $u$ satisfies \eqref{Allen-Cahn-Bratu_equation}.

The spherical domain is discretized by embedding $\Omega$ into the cube $[-1,1]^3$, which is partitioned into $2\times2\times2$ subdomains. 
Collocation points are uniformly distributed in each subdomain, and only those with $x^2+y^2+z^2<1$ are retained to assemble the NLS system. 
To impose Dirichlet conditions on $\partial\Omega$, 30,000 uniformly distributed boundary points are generated on the unit sphere using the Fibonacci spiral method.

Table~\ref{3D_PDE_results_Allen-Cahn-Bratu} reports the numerical results of the AMIPN method. 
With the subdomain partition fixed at $(N_x,N_y,N_z)=(2,2,2)$ and $\sqrt[3]{Q}=30\!-\!35$, 
increasing the number of random basis functions $J$ from 400 to 2000 leads to substantial accuracy improvement. 
Specifically, the relative $L^2$ error of $u$ decreases from $2.21\times10^{-4}$ to $1.51\times10^{-9}$, 
while the errors of its first derivatives $(u_x,u_y,u_z)$ decrease from $O(10^{-3}\!-\!10^{-4})$ to $O(10^{-8})$. 
The computational time increases moderately from about 84\,s to 1287\,s, 
whereas the numbers of outer iterations and Jacobian evaluations remain fixed at 3 and 2, respectively.

\begin{table}[!htbp]
    \centering
    \caption{Results of the AMIPN method for the 3D 3D steady-state nonlinear diffusion-reaction equation \eqref{Allen-Cahn-Bratu_equation}}
    \label{3D_PDE_results_Allen-Cahn-Bratu}
    \begin{tabular}{c c c c c c c c c c c}
      \toprule
      $(N_x, N_y, N_z)$ & $\sqrt[3]{Q}$ & $J$ & Solver & IT & NJ & CPU& $u$ error & $u_x$ error & $u_y$ error & $u_z$ error \\
      \midrule
      \multirow{5}{*}{$(2,2, 2)$} 
        & \multirow{5}{*}{$30$} & \multirow{1}{*}{$400$} & AMIPN & $3$ & $2$ &$83.54$ & $2.21e-04$ & $1.00e-03$ & $ 9.94e-04$ & $9.77e-04$ \\
        \cmidrule(lr){3-11}
        && \multirow{1}{*}{$800$}  &AMIPN & $3$ & $2$ &$193.29$  & $3.83e-06$ & $2.42e-05$ & $2.40e-05$ & $ 2.34e-05$ \\
       
        \cmidrule(lr){3-11}
        && \multirow{1}{*}{$1200$} &AMIPN & $3$ & $2$  &$348.00$ & $2.06e-07$ & $1.53e-06$ & $1.53e-06$ & $1.52e-06$ \\
       
        \cmidrule(lr){3-11}
        && \multirow{1}{*}{$1600$} &AMIPN & $3$ & $2$  &$ 616.20$ & $1.62e-08$ & $1.35e-07$ & $1.43e-07$ & $1.35e-07$ \\
        \cmidrule(lr){2-11}
          & \multirow{1}{*}{$35$}& \multirow{1}{*}{$2000$}&AMIPN & $3$ & $2$  &$1286.85$ & $1.51e-09$ & $1.41e-08$ & $1.41e-08$ & $1.42e-08$ \\
      \bottomrule
    \end{tabular}
\end{table}

At the end of this section, we consider a larger-scale benchmark problem: 
the three-dimensional steady-state Gray-Scott reaction-diffusion system posed on the cubic domain $\Omega=(0,1)^3$:
\begin{equation}\label{3D_Gray-Scott_equation}
\begin{cases}
D_u \Delta u - u v^2 + F(1-u) = f_1(x,y,z), \\[3pt]
D_v \Delta v + u v^2 - (F+k)v = f_2(x,y,z),
\end{cases}
\quad (x,y,z)\in\Omega,
\end{equation}
with the exact solution
\[
u(x,y,z)=\sin(x(1-x))\sin(y(1-y))e^z, 
\quad 
v(x,y,z)=\sin(\pi x)\sin(\pi y)\sin(\pi z).
\]
In the numerical experiments, we set $F=0.060$, $k=0.062$, and use solution-dependent diffusion coefficients 
$D_u=1+u^2$ and $D_v=1+v^2$. Homogeneous Dirichlet boundary conditions are imposed on $\partial\Omega$. 
When discretized by the RFM, this coupled system yields a nonlinear least-squares problem whose Jacobian 
has a dimension approximately twice that of the single-equation problems 
\eqref{3D_nonlinear_PDE}-\eqref{3D_nonlinear_Helmholtz_equation}.

Table~\ref{3D_PDE_results_Gray-Scott} shows that the AMIPN method maintains excellent stability 
and scalability for this large-scale strongly nonlinear coupled system. 
As the numbers of collocation points and basis functions increase, the relative errors of $u$ and $v$ 
decrease from $O(10^{-5})$ to $O(10^{-10})$, exhibiting spectral-like convergence. 
Convergence is achieved within only 3-4 outer iterations and 2 Jacobian evaluations, 
highlighting the strong potential of AMIPN for large-scale nonlinear PDE systems.

\begin{table}[!htbp]
  \centering
  \caption{Results of the AMIPN method for the 3D steady-state Gray-Scott reaction-diffusion system \eqref{3D_Gray-Scott_equation}}
  \label{3D_PDE_results_Gray-Scott}
  \begin{tabular}{c c c c c c c c}
    \toprule
    $(N_x, N_y, N_z)$ & $\sqrt[3]{Q} $ & $J$ & Solver & IT &NJ & $u$ error & $v$ error \\
    \midrule
    \multirow{7}{*}{$(2, 2, 2)$} 
        & \multirow{2}{*}{$20$} &  \multirow{1}{*}{400}& AMIPN & $3$ &$2$ &$7.60e-05$ & $ 5.00e-05$ \\
         \cmidrule(lr){3-8}
       && \multirow{1}{*}{800}  & AMIPN & $3$ &$2$ &$ 1.35e-06$ & $7.37e-07$ \\
    \cmidrule(lr){2-8}
       & \multirow{2}{*}{$25$} &  \multirow{1}{*}{1200} & AMIPN & $3$ &$2$ &$ 6.08e-08$ & $3.71e-08$ \\
         \cmidrule(lr){3-8}
       && \multirow{1}{*}{1600}  & AMIPN & $3$ &$2$ &$3.78e-09$ & $3.13e-09$ \\
      \cmidrule(lr){2-8}
      & \multirow{2}{*}{$30$} &  \multirow{1}{*}{1600} & AMIPN & $4$ &$2$ &$ 5.01e-09$ & $2.90e-09$ \\
    \cmidrule(lr){3-8}
       && \multirow{1}{*}{2000}  & AMIPN & $4$ &$2$  &$ 4.05e-10$ & $ 4.56e-10$ \\ 
    \bottomrule
  \end{tabular}
\end{table}

\subsection{Time-dependent nonlinear PDEs on complex evolving domains}
Time-dependent PDEs on complex space-time domains pose significant challenges, since they often involve moving boundaries, topological changes, large deformations, and complex geometries. In such scenarios, conventional mesh-based methods typically suffer from severe mesh distortion and frequent remeshing, leading to substantial implementation difficulties, loss of accuracy, and high computational cost.

Owing to its mesh-free nature, the proposed approach can accurately handle evolving geometries 
while maintaining high accuracy and numerical stability. 
To systematically evaluate its performance, we consider four benchmark problems with progressively increasing geometric complexity:  
(i) the two-dimensional Allen-Cahn equation with a moving interior hole;  
(ii) the two-dimensional Klein-Gordon equation with moving internal boundaries and topological changes;  
(iii) a nonlinear reaction-diffusion-convection system defined on a time-dependent domain undergoing large deformations; and  
(iv) a nonlinear diffusion equation and a Lotka-Volterra reaction-diffusion system, defined on separate complex space-time domains, thereby providing complementary benchmarks.

We first consider the Allen-Cahn equation defined on the three-dimensional space-time domain $\left([0,1]^2 \setminus \mathcal{H}(t)\right)\times(0,1]$:
\begin{equation}\label{eq:allen_cahn_movinghole}
\begin{cases}
\frac{\partial u}{\partial t} - \Delta u + u^3 - u = f(x,y,t), & (x,y,t)\in([0,1]^2\setminus\mathcal{H}(t))\times(0,1], \\[3pt]
u(x,y,0)=u_0(x,y), & (x,y)\in[0,1]^2\setminus\mathcal{H}(0), \\[3pt]
u(x,y,t)=g(x,y,t), & (x,y,t)\in\partial[0,1]^2\cup\partial\mathcal{H}(t),~t\in[0,1],
\end{cases}
\end{equation}
where $\Delta u=u_{xx}+u_{yy}$ and $\mathcal{H}(t)=\{(x,y):(x-c_x(t))^2+(y-c_y(t))^2\le r^2\}$ 
denotes a circular hole of radius $r=0.1$ with center
$c_x(t)=0.5+0.2\cos(\pi t)$, $c_y(t)=0.5+0.2\sin(\pi t)$. 
The source term $f$ is chosen so that the exact solution is
\[
u(x,y,t)=2\sin(x(1-x))\sin(y(1-y))e^{t+1}.
\]

\begin{figure}[!htbp]
    \renewcommand\figurename{Figure}
      \centering
      \includegraphics[scale=0.35]{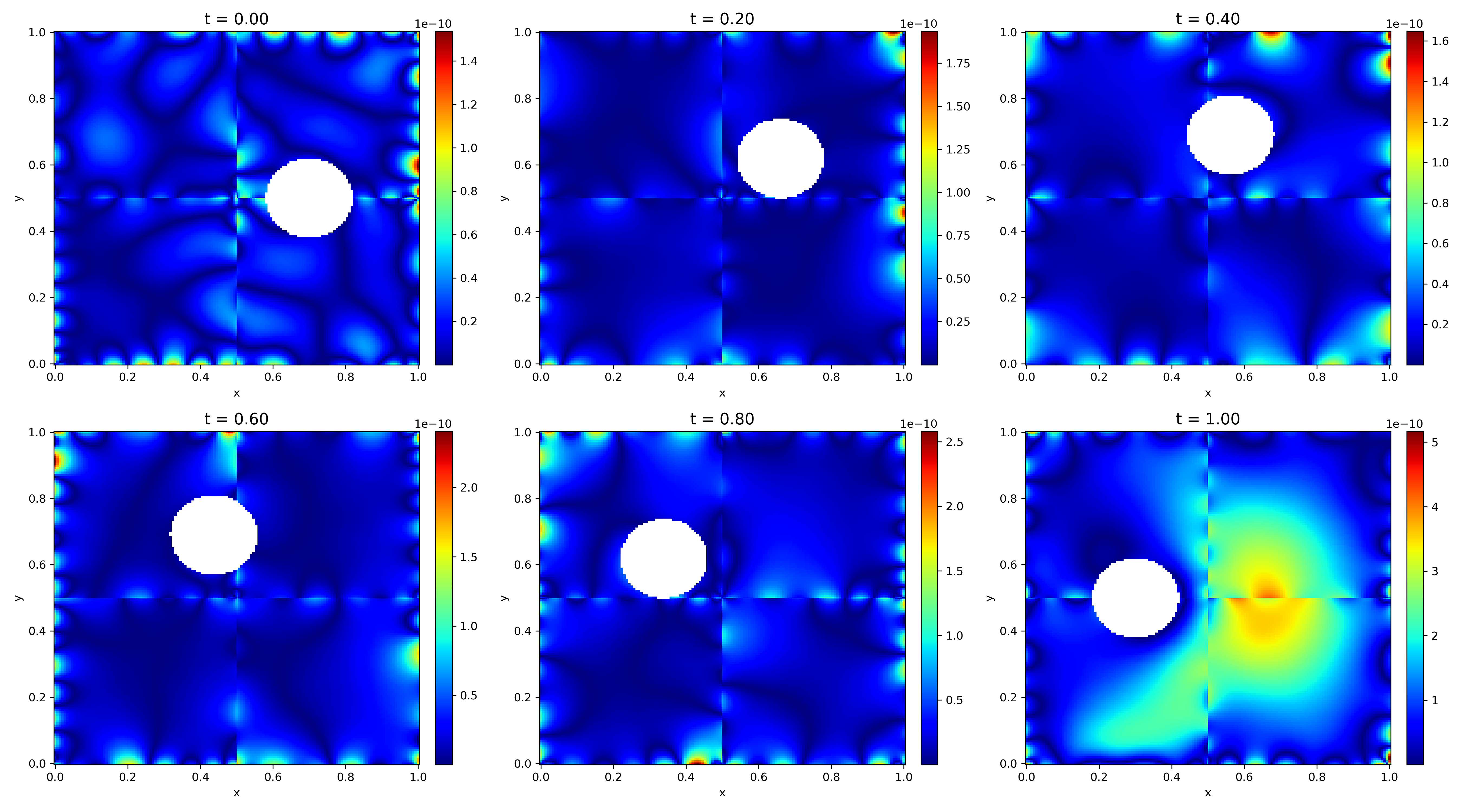}
\caption{Absolute error on the $x$-$y$ plane at times $t = 0, 0.20, 0.40, 0.60, 0.80, \text{and} 1.00$ for the Allen--Cahn equation \eqref{eq:allen_cahn_movinghole}}\label{fig:Allen-Cahn_surface_error}
\end{figure}

Table~\ref{2D_Allen-Cahn} reports the numerical results of the AMIPN method for this moving-hole problem. 
As the numbers of collocation points $\sqrt[3]{Q}$ per subdomain and random basis functions $J$ increase, 
the $L^2$ error of $u$ decreases from $O(10^{-5})$ to $O(10^{-10})$, 
while the errors of its derivatives $(u_x,u_y,u_t)$ decrease from $O(10^{-3}\!-\!10^{-4})$ to $O(10^{-8}\!-\!10^{-10})$, 
exhibiting a nearly exponential convergence trend. 
The number of outer iterations remains fixed at three, with one Jacobian evaluation per iteration, 
demonstrating the efficiency and stability of the method. 
The computational time increases predictably with $J$, from about 23\,s to 1066\,s. 
Figure~\ref{fig:Allen-Cahn_surface_error} visualizes the absolute errors at $t=0,0.20,0.40,0.60,0.80$, and $1.00$, 
clearly illustrating both the temporal evolution of the circular hole and the accuracy of the method.

\begin{table}[!htbp]
  \centering
  \caption{Results of the AMIPN method for the 2D Allen-Cahn equation \eqref{eq:allen_cahn_movinghole}}
  \label{2D_Allen-Cahn}
  \begin{tabular}{c c c c c c c c c c c}
    \toprule
    $(N_x, N_y, N_z)$ & $\sqrt[3]{Q} $ & $J$ & Solver & IT& NJ & CPU & $u$ error & $u_x$ error & $u_y$ error & $u_t$ error\\
    \midrule
    \multirow{7}{*}{$(2, 2, 2)$} 
        & \multirow{2}{*}{$20$} & \multirow{1}{*}{400} & AMIPN & $3$ & $1$ & $22.78$ & $7.70e-05$ & $3.13e-04$ & $3.18e-04$& $3.26e-03$\\
        && \multirow{1}{*}{800} & AMIPN & $3$ & $1$ & $56.50$ & $1.73e-06$& $ 8.78e-06$ & $8.62e-06$& $ 1.11e-04$\\
    \cmidrule(lr){2-11}
        & \multirow{2}{*}{$25$} 
        &\multirow{1}{*}{1200} & AMIPN & $3$ & $1$ & $ 159.11$ & $ 5.47e-08$ & $3.59e-07$ & $3.54e-07$& $4.04e-06$\\
        && \multirow{1}{*}{1600} & AMIPN & $3$ & $1$ & $282.21$ & $ 5.30e-09$ & $3.25e-08$ & $3.23e-08$& $ 5.35e-07$\\
      \cmidrule(lr){2-11}
        & \multirow{3}{*}{$30$} 
        & \multirow{1}{*}{1600} & AMIPN & $3$ & $1$ & $393.06$ & $4.01e-09$& $2.59e-08$ & $2.67e-08$& $ 3.60e-07$\\
        && \multirow{1}{*}{2000} & AMIPN & $3$ & $1$ & $623.36$ & $ 7.25e-10$& $5.22e-09$ & $5.30e-09$& $7.77e-08$\\
        && \multirow{1}{*}{2400} & AMIPN & $3$ & $1$ & $ 1066.48$ & $1.82e-10$& $8.16e-10$ & $7.92e-10$& $ 1.33e-08$\\
    \bottomrule
  \end{tabular}
\end{table}

\begin{figure}[!htbp]
    \renewcommand\figurename{Figure}
      \centering
      \includegraphics[scale=0.35]{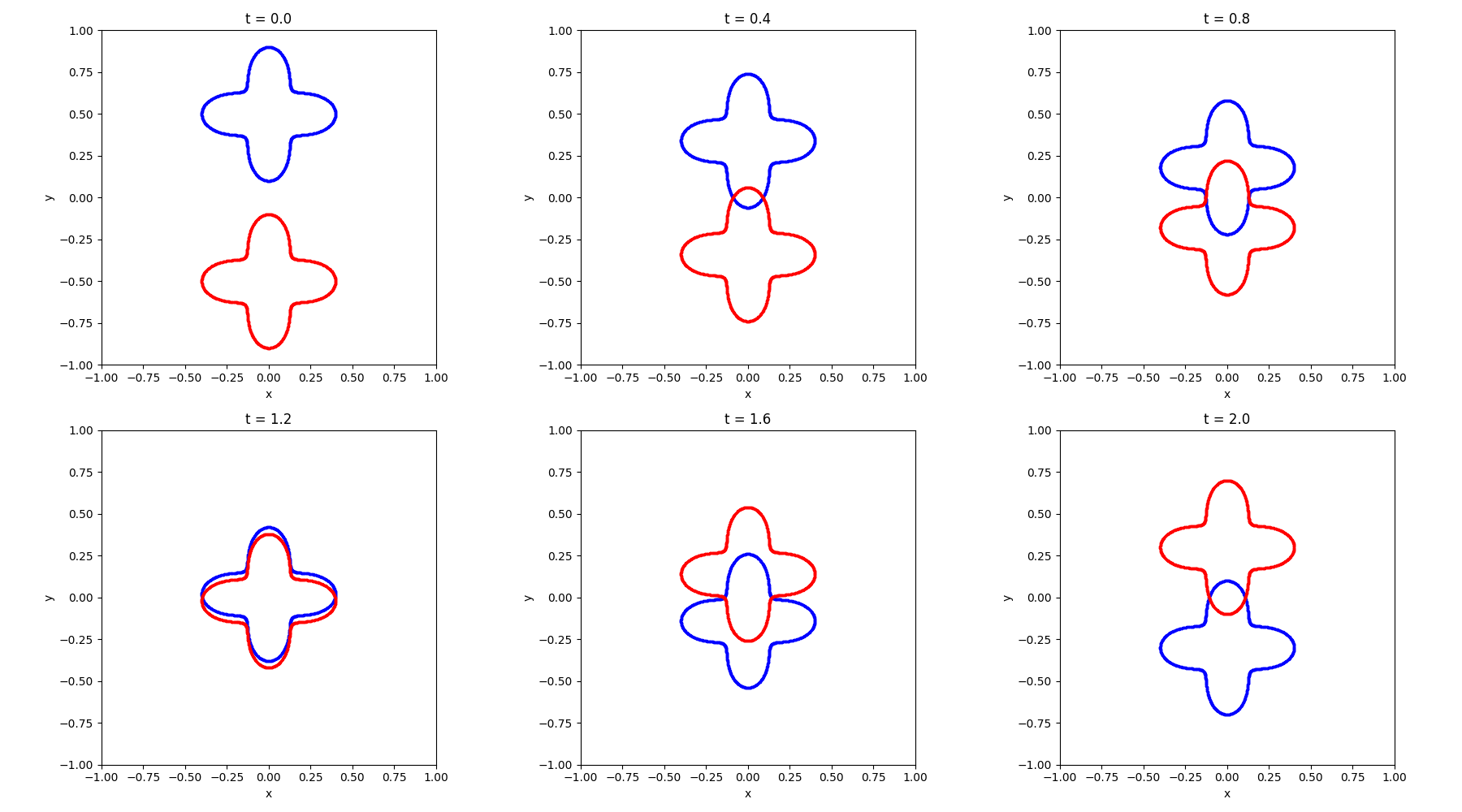}
    \caption{The evolution of the zero level sets of $F_1(x,y,t)$ and $F_2(x,y,t)$, where the corresponding excised regions first merge and then split over time.}
    \label{fig:merging_process}
\end{figure}

We next consider the two-dimensional Klein-Gordon equation posed on a square spatial domain containing two dynamically evolving internal obstacles:
\begin{equation}\label{eq:KG_system}
\left\{
\begin{aligned}
&u_{tt}-\Delta u + u + \beta u^2 = f(x,y,t),
&& (x,y,t)\in \big([-1,1]^2 \setminus (\Omega_1(t)\cup\Omega_2(t))\big)\times (0,2], \\[3pt]
&u(x,y,0)=\phi(x,y), \quad u_t(x,y,0)=\psi(x,y),
&& (x,y)\in [-1,1]^2\setminus(\Omega_1(0)\cup\Omega_2(0)), \\[3pt]
&u(x,y,t)=h(x,y,t),
&& (x,y,t)\in \partial([-1,1]^2)\cup \partial\Omega_1(t)\cup\partial\Omega_2(t),\; t\in[0,2],
\end{aligned}
\right.
\end{equation}
where $\Delta u = u_{xx}+u_{yy}$, and $\beta$ is a prescribed nonlinear coefficient. 
Each internal obstacle is defined by a level set $\Omega_i(t)=\{(x,y)\mid F_i(x,y,t)<0\}$ ($i=1,2$). 
During the simulation, the two obstacles move toward each other, merge into a single connected component, and later separate again, thereby inducing a topological change in the computational domain. 
The level-set functions are
\begin{equation}\label{eq:flower_geometries}
\begin{aligned}
F_1(x,y,t) &= x^2 + \big(y - 0.5 + 0.4t\big)^2 - \Big(0.3 + 0.1 \cos\big(4\theta_1\big)\Big)^2, \\
F_2(x,y,t) &= x^2 + \big(y + 0.5 - 0.4t\big)^2 - \Big(0.3 + 0.1 \cos\big(4\theta_2\big)\Big)^2, \\
\theta_1 &= \arctan\!\Big(\frac{y - 0.5 + 0.4t}{x}\Big),
\qquad
\theta_2 = \arctan\!\Big(\frac{y + 0.5 - 0.4t}{x}\Big).
\end{aligned}
\end{equation}
Figure~\ref{fig:merging_process} visualizes the merging process at $t = 0,\,0.40,\,0.80,\,1.20,\,1.60,\,\text{and}\,2.00$. 
For verification, the source term $f$, initial data $\phi,\psi$, and boundary data $h$ are chosen consistently. 
The exact solution is given by
\begin{equation}\label{eq:KG_exact}
u(x,y,t)=\frac{x+y+t+2}{1+x^2+y^2}.
\end{equation}

\begin{figure}[!htbp]
    \renewcommand\figurename{Figure}
      \centering
      \includegraphics[scale=0.35]{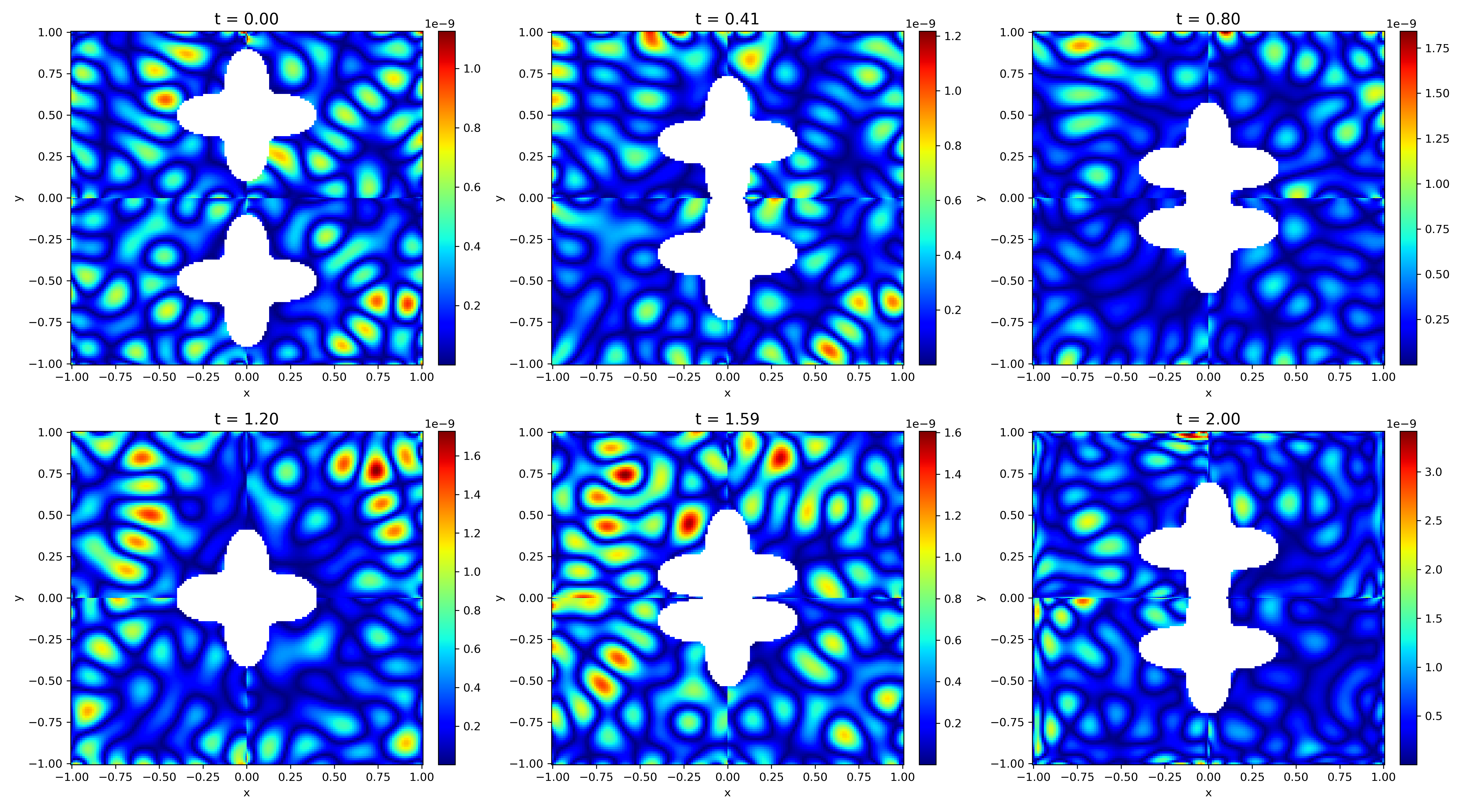}
    \caption{Absolute error on the spatial cross-sections at $t = 0$, $0.41$, $0.80$, $1.20$, $1.59$, and $2.00$ for the Klein-Gordon equation \eqref{eq:KG_system}}
    \label{fig:Klein-Gordon_surface_error}
\end{figure}

Table~\ref{Klein-Gordon_result} summarizes the numerical results of the AMIPN method. 
With the subdomain partition fixed at $(N_x,N_y,N_z)=(2,2,2)$ and collocation density $\sqrt[3]{Q}=30$, 
increasing the number of random features $J$ from 400 to 2400 leads to substantial accuracy improvement. 
The relative error of $u$ decreases from $2.93\times10^{-5}$ to $2.12\times10^{-10}$, 
while the errors in $u_x$ and $u_y$ decrease from $O(10^{-4})$ to $O(10^{-9})$, 
and the $u_t$ error decreases from $O(10^{-3}\!-\!10^{-4})$ to $O(10^{-8})$, 
exhibiting spectral-like convergence. 
The number of outer iterations remains fixed at four, with three Jacobian evaluations per iteration, 
while the computational time increases from about 123\,s to 2101\,s as $J$ increases. 
Figure~\ref{fig:Klein-Gordon_surface_error} visualizes the absolute errors of the AMIPN method at different time instances, 
highlighting its robustness for problems with evolving internal obstacles and topological changes, 
and underscoring its capability to handle complex space-time domains with dynamic topology.

\begin{table}[!htbp]
  \centering
  \caption{The results of the AMIPN method for the Klein-Gordon equation \eqref{eq:KG_system}}
  \label{Klein-Gordon_result}
  \begin{tabular}{c c c c c c c c c c c}
    \toprule
    $(N_x, N_y, N_z)$ & $\sqrt[3]{Q} $ & $J$ & Solver & IT & NJ & CPU & $u$ error & $u_x$ error& $u_y$ error& $u_t$ error\\
    \midrule
    \multirow{7}{*}{$(2, 2, 2)$} 
        & \multirow{7}{*}{$30$} & \multirow{1}{*}{400} & AMIPN& $4$ & $3$& $123.12$ & $ 2.93e-05$ & $3.52e-04$& $3.84e-04$& $9.25e-04$\\
    \cmidrule(lr){4-11}
        & 
        &\multirow{1}{*}{800} & AMIPN & $4$ & $3$& $330.70$ & $8.93e-07$ & $1.34e-05$& $ 1.47e-05$& $3.67e-05$\\
        \cmidrule(lr){4-11}
         && \multirow{1}{*}{1200} & AMIPN & $4$ & $3$& $582.42$ & $ 6.67e-08$ & $1.27e-06$& $ 1.28e-06$& $ 3.38e-06$\\
      \cmidrule(lr){4-11}
        & 
        & \multirow{1}{*}{1600} & AMIPN & $4$ & $3$& $890.34$ & $7.56e-09$ & $ 1.41e-07$& $1.71e-07$& $4.28e-07$\\
         \cmidrule(lr){4-11}
        && \multirow{1}{*}{2000} & AMIPN & $4$ & $3$& $ 1496.10$ & $ 1.02e-09$ & $ 2.18e-08$& $2.36e-08$& $6.17e-08$\\
        \cmidrule(lr){4-11}
        && \multirow{1}{*}{2400} & AMIPN & $4$ & $3$& $2100.97$ & $2.12e-10$& $4.90e-09$& $5.40e-09$& $ 1.42e-08$\\
    \bottomrule
  \end{tabular}
\end{table}

We next consider a nonlinear reaction-diffusion-convection system defined on a time-dependent perforated domain with a deforming interior hole.
The governing equations are given by
\begin{equation}\label{eq:rdc_large_def}
\left\{
\begin{aligned}
&u_t+\mathbf{w}(x,y,t)\!\cdot\!\nabla u - \Delta u + \alpha\,u v = f_1(x,y,t),
&& (x,y,t)\in\big((0,1)^2\setminus\Omega(t)\big)\times(0,1],\\[3pt]
&v_t+\mathbf{w}(x,y,t)\!\cdot\!\nabla v - \Delta v + \beta\,e^{u} v^{2} = f_2(x,y,t),
&& (x,y,t)\in\big((0,1)^2\setminus\Omega(t)\big)\times(0,1],\\[3pt]
&u(x,y,0)=\phi(x,y),\quad v(x,y,0)=\psi(x,y),
&& (x,y)\in(0,1)^2\setminus\Omega(0),\\[3pt]
&u(x,y,t)=h_1(x,y,t),\quad v(x,y,t)=h_2(x,y,t),
&& (x,y,t)\in\partial(0,1)^2\cup\partial\Omega(t),\ t\in[0,1].
\end{aligned}
\right.
\end{equation}
Here $\Delta=\partial_{xx}+\partial_{yy}$ and we set $\alpha=\beta=1$. 
The prescribed advection velocity field is
\begin{align*}
\mathbf{w}(x, y, t)=\cos \left(\frac{\pi t}{10}\right)\left(\sin ^2(\pi x) \sin (2 \pi y),-\sin ^2(\pi y) \sin (2 \pi x)\right)^{\top}.
\end{align*}
The computational domain is the unit square $(0,1)^2$ with a time-dependent interior hole $\Omega(t)$. 
It is initially defined as
\begin{align*}
\Omega(0)=\left\{(x, y)\in(0,1)^2 \mid(x-0.5)^2+(y-0.75)^2 \leq 0.15^2\right\}.
\end{align*}
The interior perforation $\Omega(t)$ evolves over time according to the prescribed advection velocity field $\mathbf{w}(x, y, t)$. 
To track the motion of the perforation boundary, the boundary points $\mathbf{X}(t)\in\partial\Omega(t)$ are advected by the flow field according to
\[
\frac{d \mathbf{X}}{d t}=\mathbf{w}(\mathbf{X}(t), t),\qquad \mathbf{X}(0)\in\partial\Omega(0).
\]
The trajectories of boundary points are numerically integrated using a fifth-order Runge-Kutta scheme. 
At each selected time $t$, the updated positions of the advected boundary points are used to reconstruct the moving interface $\partial \Omega(t)$, thereby defining the time-dependent perforated computational domain $(0,1)^2\setminus\Omega(t)$.
For verification, the source terms $f_1,f_2$ and the boundary/initial data are chosen consistently with the exact solution
\begin{align*}
u(x, y, t)=\sin (x(1-x)) \sin (y(1-y)) e^t, \qquad 
v(x, y, t)=\sin (\pi x) \sin (\pi y) \sin (\pi t).
\end{align*}

Table~\ref{2D_large_deformation} reports the numerical results of the AMIPN method for the nonlinear reaction-diffusion-convection system on a time-dependent domain with a deforming interior hole. 
The relative $L^2$ errors of $u$ and $v$ decrease from $O(10^{-5}\!-\!10^{-6})$ to $O(10^{-10})$, 
exhibiting spectral-like convergence and confirming the spectral accuracy of the proposed method for this challenging large-deformation problem. 
Moreover, even in large-scale and potentially ill-conditioned settings, 
the number of outer iterations remains fixed at three, with only two Jacobian evaluations per iteration, 
demonstrating the efficiency and robustness of the proposed approach.

\begin{figure}[!htbp]
    \renewcommand\figurename{Figure}
    \centering
    \begin{subfigure}[t]{0.48\textwidth}
        \centering
        \includegraphics[width=\textwidth]{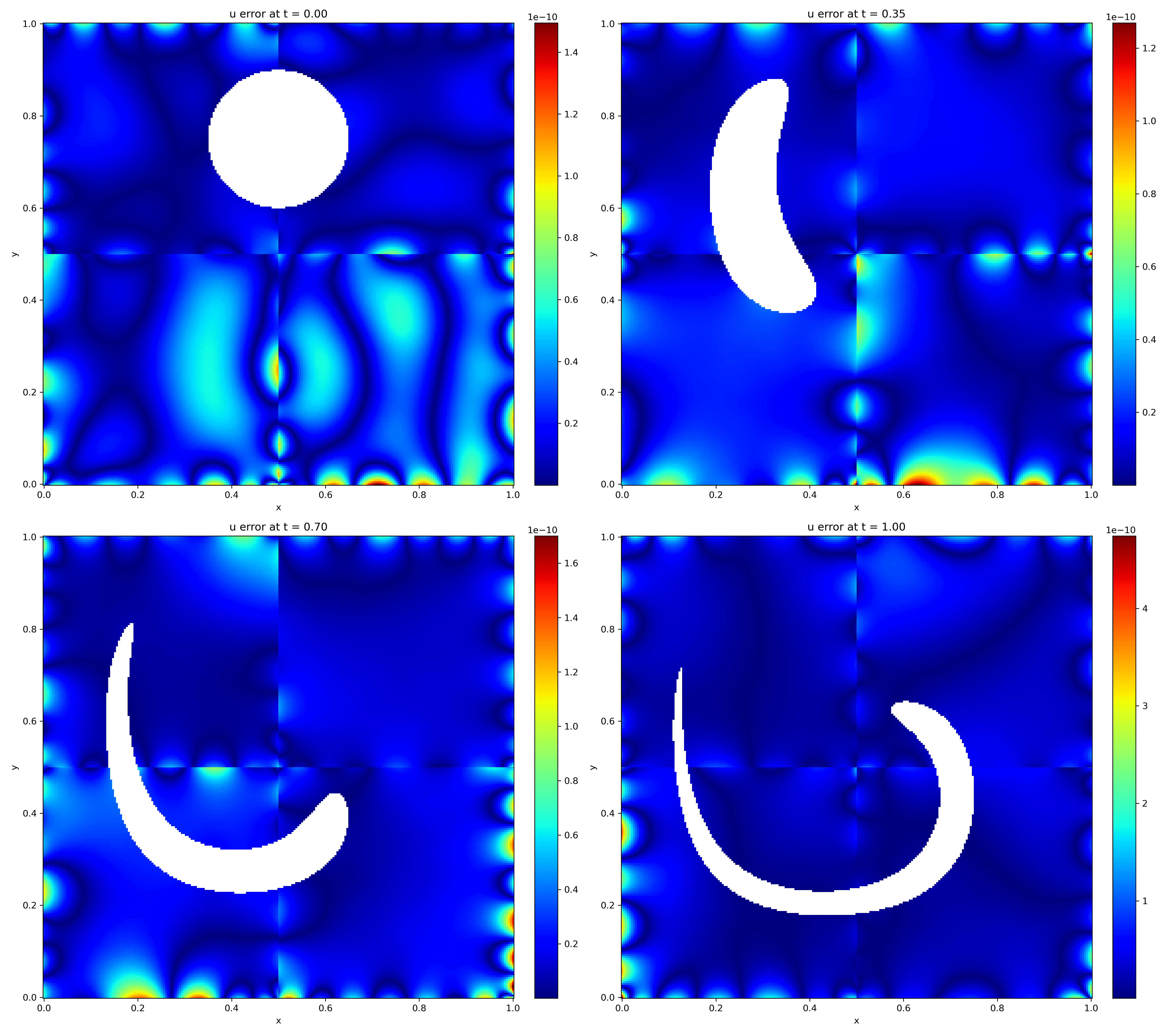}
        \caption{$u$ error}
        \label{fig:u_error}
    \end{subfigure}
    \hfill
    \begin{subfigure}[t]{0.48\textwidth}
        \centering
        \includegraphics[width=\textwidth]{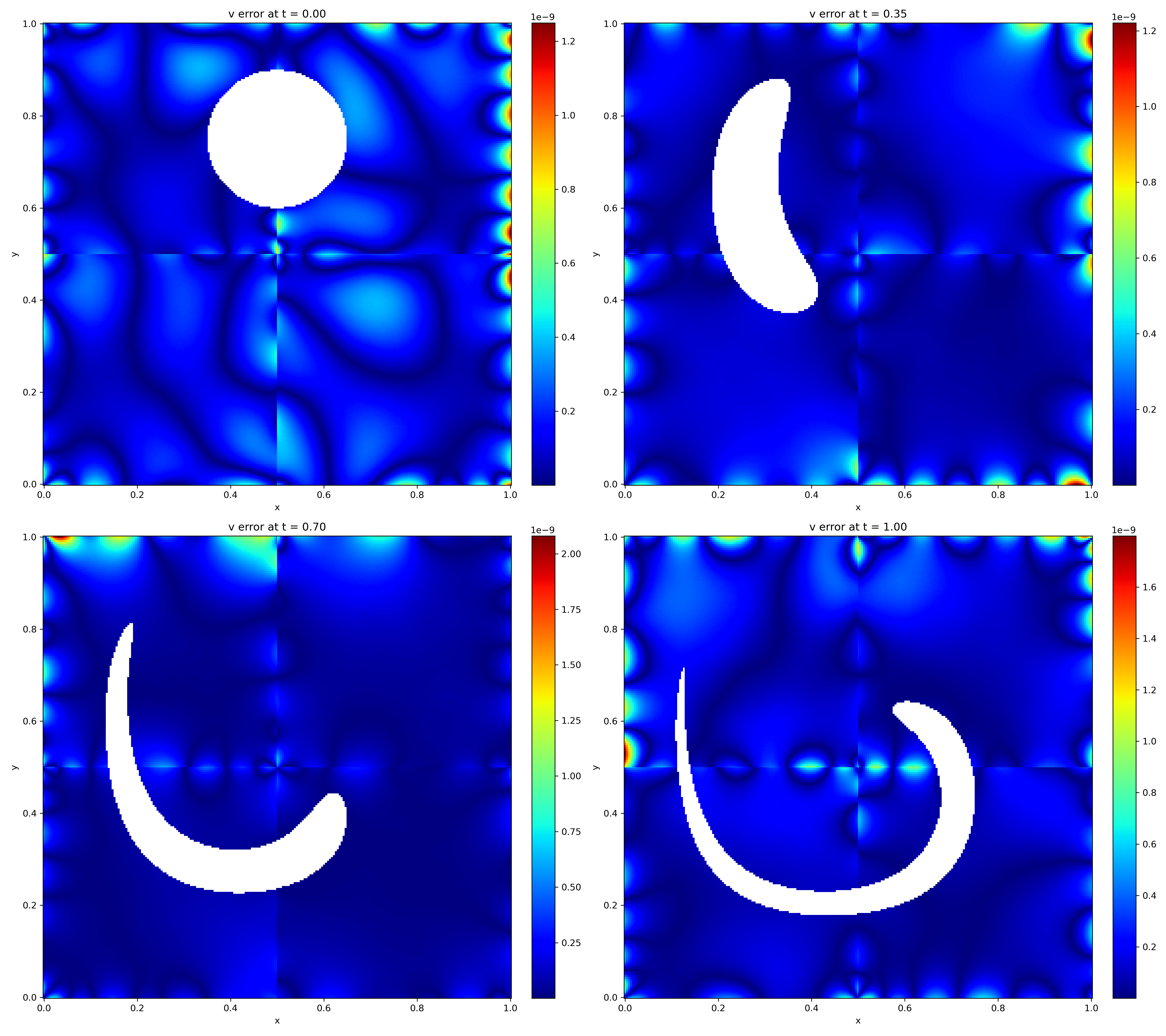}
        \caption{$v$ error}
        \label{fig:v_error}
    \end{subfigure}
    \caption{Absolute errors of $u$ and $v$ at different time instances $t=0,~0.35,~0.70,~1.00$. The white region indicates the moving hole $\Omega(t)$ transported by the prescribed velocity field.}
    \label{fig:uv_error_surface}
\end{figure}

Table~\ref{PCPU_LSQR_Iter} reports the preconditioning time per outer iteration and the number of LSQR iterations in the inner loop. 
The results indicate that the preconditioning stage constitutes the dominant portion of the computational cost. 
By reusing the preconditioned Jacobian across multiple inner iterations, 
the AMIPN method substantially reduces the number of preconditionings, 
thereby reducing the overall cost compared with IPN. 
Figure~\ref{fig:uv_error_surface} illustrates the absolute errors of $u$ and $v$ at $t=0,0.35,0.70,1.00$. 
The white region denotes the interior hole $\Omega(t)$, which deforms from an initial circular shape into a highly elongated snake-like geometry under the prescribed advection field, 
clearly demonstrating the capability of AMIPN to maintain high accuracy and stability in strongly deforming domains with time-dependent geometries.

\begin{table}[!htbp]
  \centering
  \caption{Results of the AMIPN method for the reaction-diffusion-convection system \eqref{eq:rdc_large_def}}
  \label{2D_large_deformation}
  \begin{tabular}{c c c c c c c c c}
    \toprule
    $(N_x, N_y, N_z)$ & $\sqrt[3]{Q} $ & $J$ & Solver & IT &NJ & CPU & $u$ error & $v$ error \\
    \midrule
    \multirow{7}{*}{$(2, 2, 2)$} 
        & \multirow{2}{*}{$20$} &  \multirow{1}{*}{400} & AMIPN & $3$ &$2$ &$141.81$ &$6.99e-05$ & $ 5.48e-05$ \\
       && \multirow{1}{*}{800}  & AMIPN & $3$ &$2$ &$422.22$ &$ 1.53e-06$ & $ 8.28e-07$ \\
    \cmidrule(lr){2-9}
       & \multirow{2}{*}{$25$} &  \multirow{1}{*}{1200} & AMIPN & $3$ &$2$ &$ 1294.65$ &$ 6.04e-08$ & $3.91e-08$ \\
       && \multirow{1}{*}{1600}  & AMIPN & $3$ &$2$ &$ 2817.73$ &$3.88e-09$ & $3.99e-09$ \\
      \cmidrule(lr){2-9}
      & \multirow{2}{*}{$30$} &  \multirow{1}{*}{1600} & AMIPN & $3$ &$2$ &$42 44.36$ &$5.03e-09$ & $3.25e-09$ \\
       && \multirow{1}{*}{2000}  & AMIPN & $3$ &$2$ &$5438.12$ &$4.76e-10$ & $5.48e-10$ \\
    \bottomrule
  \end{tabular}
\end{table}
\begin{table}[!htbp]
  \centering
  \caption{Preconditioning time per outer iteration and the number of LSQR steps in the inner iteration}
  \label{PCPU_LSQR_Iter}
  \begin{tabular}{c c c c c c c c c}
    \toprule
    $(N_x, N_y, N_z)$ & $\sqrt[3]{Q} $ & $J$ &\textbf{Outer Iter} & \textbf{PCPU$(s)$} & \textbf{Iter 1} & \textbf{Iter 2} & \textbf{Iter 3} & \textbf{TLI}\\
    \midrule
    \multirow{6}{*}{$(2, 2, 2)$} 
        & \multirow{6}{*}{$25$} &  \multirow{3}{*}{1200}& $1$ & $ 476.95$ & $33$ & $31$ & $34$ & $98$ \\
       &&& $2$ & $ 478.64$ & $32$ & $11$ & $--$ & $43$ \\
       &&& $3$ &$--$ & $9$ & $--$ & $--$ & $9$ \\
         \cmidrule(lr){3-9}
       && \multirow{3}{*}{1600}&$1$ & $1072.94$ & $35$ & $34$ & $35$ & $104$ \\
       &&& $2$ & $1052.04$ & $39$ & $16$ & $--$ & $55$ \\
       &&& $3$ &$--$ & $14$ & $--$ & $--$ & $14$ \\
    \bottomrule
  \end{tabular}
\end{table}

We next consider a Lotka-Volterra reaction-diffusion system defined on a perforated space-time domain:
\begin{equation}\label{eq:lv_complex}
\begin{cases}
u_t - \Delta u -u + u v = f_1(x,y,t), & (x,y,t)\in \mathcal{D}^*,\\
v_t - \Delta v + v - u v = f_2(x,y,t), & (x,y,t)\in \mathcal{D}^*,\\
u(x,y,0)=u_0(x,y),\quad v(x,y,0)=v_0(x,y), & (x,y)\in\Omega,\\
u(x,y,t)=g_1(x,y,t),\quad v(x,y,t)=g_2(x,y,t), & (x,y,t)\in\partial\Omega\times(0,2],
\end{cases}
\end{equation}
where $\Delta u=u_{xx}+u_{yy}$ and $\Delta v=v_{xx}+v_{yy}$.
The source terms $f_1,f_2$, together with the boundary and initial data, are derived from the manufactured exact solution
\begin{align*}
    u(x, y, t)=\sin (\pi x) \sin (\pi y) e^t, \quad v(x, y, t)=\frac{x^2+y^2+t^2+2}{t+1} .
\end{align*}
The spatial background domain is the square
$
\Omega_{\mathrm{sq}}=[-1,1]^2,
$
from which a star-shaped region
\begin{align*}
    \mathcal{S}=\left\{(x, y):\left(x-c_x\right)^2+\left(y-c_y\right)^2 \leq r(\theta)^2, \quad r(\theta)=0.4+0.2 \sin (20 \theta)\right\},
\end{align*}
centered at $(c_x,c_y)=(0.02\sqrt{5},0.02\sqrt{5})$, is excised.
The resulting computational domain is
$
\Omega=\Omega_{\mathrm{sq}} \backslash \mathcal{S} .
$
Extending in time, the full space-time domain is defined as
\begin{align*}
    \mathcal{D}^*=\Omega \times(0,2]=\left([-1,1]^2 \times(0,2]\right) \backslash(\mathcal{S} \times(0,2]),
\end{align*}
which can be interpreted as a three-dimensional column $[-1,1]^2\times(0,2]$ with a star-shaped cylindrical void removed along the time axis. This construction results in a highly complex perforated space-time geometry, with oscillatory internal boundaries and nontrivial embeddings, thereby posing substantial challenges for generating high-quality meshes.

\begin{figure}[!htbp]
    \renewcommand\figurename{Figure}
    \centering
    \begin{subfigure}[t]{0.48\textwidth}
        \centering
        \includegraphics[width=\textwidth]{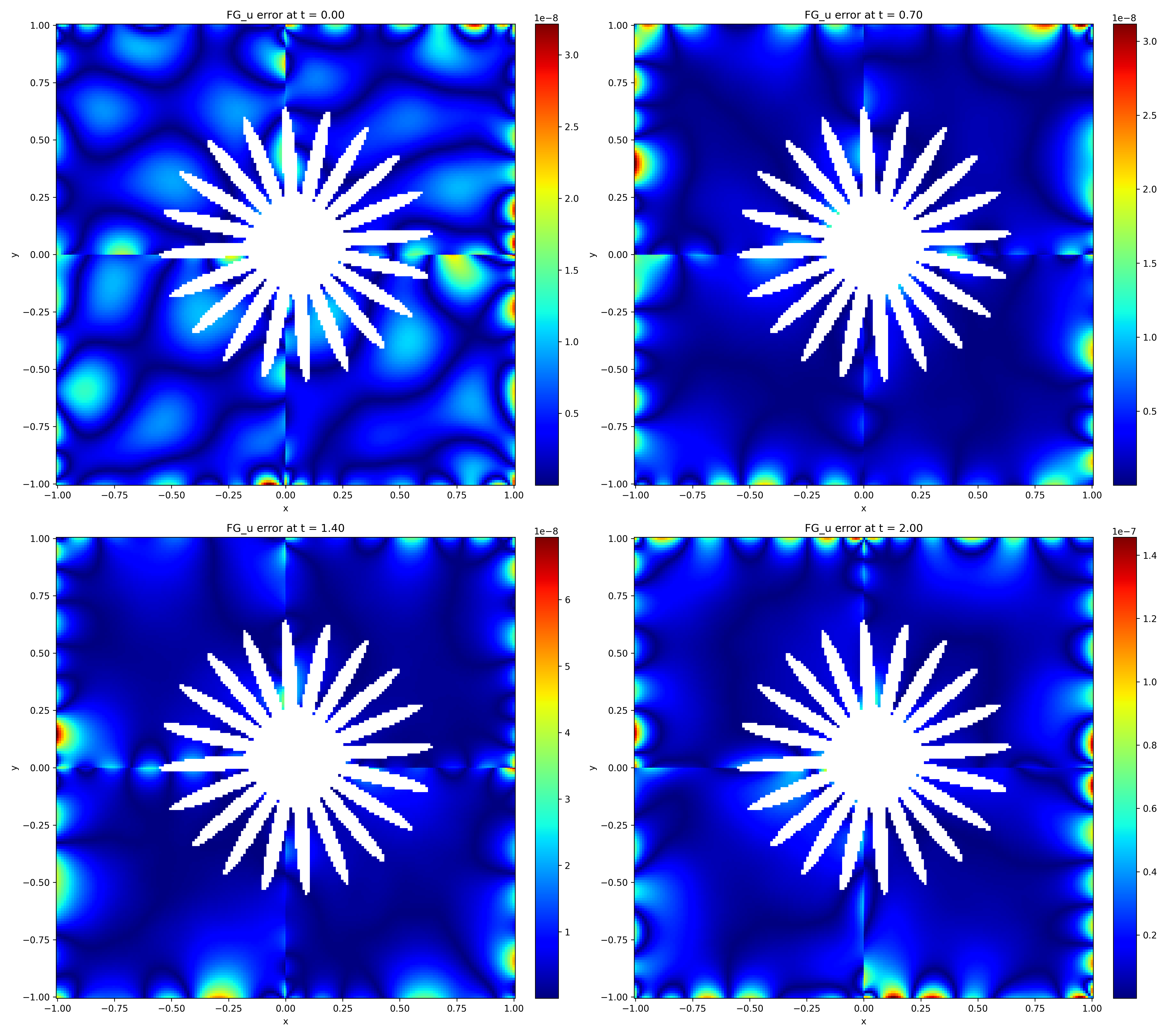}
        \caption{$u$ error}
        \label{fig:FG_u_error}
    \end{subfigure}
    \hfill
    \begin{subfigure}[t]{0.48\textwidth}
        \centering
        \includegraphics[width=\textwidth]{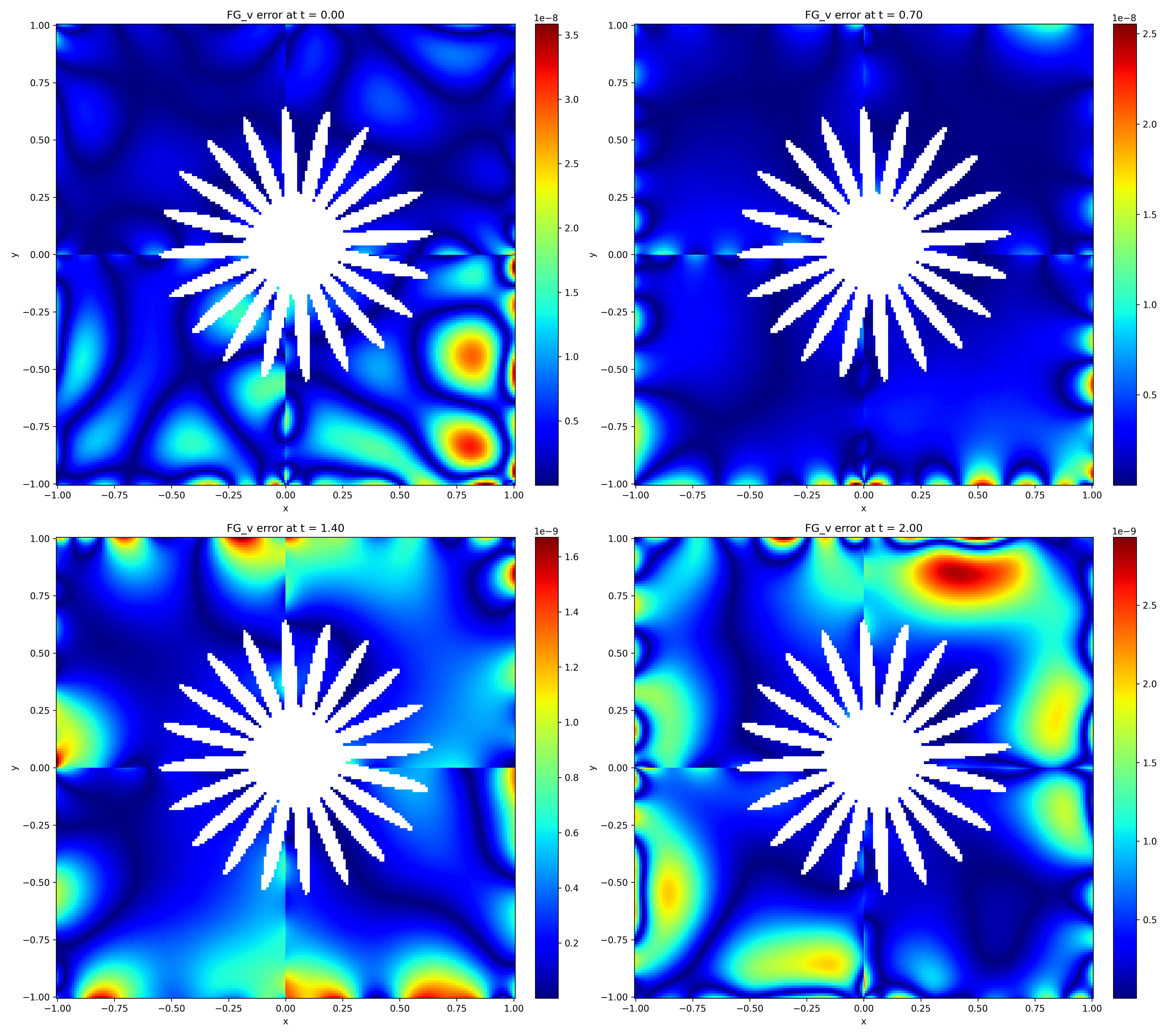}
        \caption{$v$ error}
        \label{fig:FG_v_error}
    \end{subfigure}
    \caption{Absolute errors of $u$ and $v$ at different time instances $t=0,~0.49,~1.00,~2.00$.}
    \label{fig:uv_FG_error}
\end{figure}

The results in Table \ref{2D_FG} demonstrate that the proposed AMIPN method achieves both high accuracy and stable convergence for the Lotka-Volterra reaction-diffusion system. As the number of random basis functions $J$ increases, the errors in $u$ and $v$ decrease significantly, from the order of $10^{-4}- 10^{-5}$ down to $10^{-8}-10^{-9}$. Throughout all tests, the number of outer iterations and Jacobian evaluations remain nearly constant (about 3 and 2, respectively), indicating that the computational cost does not grow with problem size. Figure \ref{fig:uv_FG_error} presents the cross-sectional errors of $u$ and $v$ at $t=0,0.7,1.40$, and $2.00$.

\begin{table}[!htbp]
  \centering
  \caption{Results of the AMIPN method for Lotka-Volterra reaction-diffusion system \eqref{eq:lv_complex}}
  \label{2D_FG}
  \begin{tabular}{c c c c c c c c}
    \toprule
    $(N_x, N_y, N_z)$ & $\sqrt[3]{Q} $ & $J$ & Solver & IT &NJ & $u$ error & $v$ error \\
    \midrule
    \multirow{7}{*}{$(2, 2, 2)$} 
        & \multirow{2}{*}{$20$} &  \multirow{1}{*}{400} & AMIPN & $3$ &$2$  &$ 3.36e-04$ & $ 3.70e-05$ \\
       && \multirow{1}{*}{800}  & AMIPN & $3$ &$2$  &$ 9.18e-06$ & $8.46e-07$ \\
    \cmidrule(lr){2-8}
       & \multirow{2}{*}{$25$} &  \multirow{1}{*}{1200} & AMIPN & $3$ &$2$ &$4.47e-07$ & $5.69e-08$ \\
       && \multirow{1}{*}{1600}  & AMIPN & $3$ &$2$ &$3.59e-08$ & $8.07e-09$ \\
      \cmidrule(lr){2-8}
      & \multirow{2}{*}{$30$} &  \multirow{1}{*}{1600} & AMIPN & $3$ &$2$ &$4.77e-08$ & $7.12e-09$ \\
       && \multirow{1}{*}{2000}  & AMIPN & $3$ &$2$ &$4.54e-09$ & $ 1.64e-09$ \\
    \bottomrule
  \end{tabular}
\end{table}

\begin{figure}[!htbp]
    \renewcommand\figurename{Figure}
        \centering
        \includegraphics[scale=0.55]{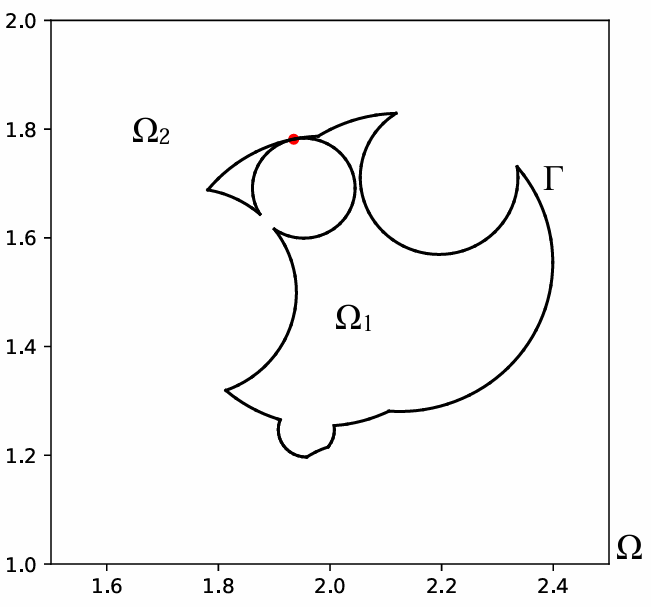}
    \caption{A two-dimensional complex geometry.}
      \label{fig:complex_geometric_COMSOL_failure}
\end{figure}

We finally consider a nonlinear diffusion problem defined on a complex space-time domain:
\begin{equation}\label{eq:complex_diffusion}
\begin{cases}
\frac{\partial u}{\partial t}
-\nabla\!\cdot\!\big((1+u^2)\nabla u\big)=f(x,y,t), & (x,y,t)\in \mathcal{D}^\ast,\\[3pt]
u(x,y,0)=u_0(x,y), & (x,y)\in\Omega,\\[3pt]
u(x,y,t)=g(x,y,t), & (x,y,t)\in\partial\Omega\times(0,1],
\end{cases}
\end{equation}
where $\Delta u=u_{xx}+u_{yy}$ and $a(u)=1+u^2$ is the nonlinear diffusion coefficient. 
The source term $f$ and the boundary/initial data $g,u_0$ are derived from the manufactured exact solution
\[
u(x,y,t)=\frac{2(x+y+t)}{1+x^2+y^2}
\Big[\tanh\!\big(\sin(\pi x)\cos(\pi y)e^{-t}\big)+1\Big].
\]
The spatial domain $\Omega$ is constructed from the background square 
$\Omega_{\mathrm{sq}}=(1.5,2.5)\times(1.0,2.0)$
by removing three circular holes 
$\mathcal{C}^-=\{B(\mathbf{c}_k,r_k)\}_{k=1}^3$
and embedding four circular inclusions 
$\mathcal{C}^+=\{B(\tilde{\mathbf{c}}_j,\tilde r_j)\}_{j=1}^4$. 
The resulting computational domain is
\[
\Omega=\big(\Omega_{\mathrm{sq}}\setminus\bigcup_{k=1}^{3}B(\mathbf{c}_k,r_k)\big)\cup\bigcup_{j=1}^{4}B(\tilde{\mathbf{c}}_j,\tilde r_j).
\]
The full space-time domain $\mathcal{D}^\ast=\Omega\times(0,1]$ may thus be interpreted as a three-dimensional column $\Omega_{\mathrm{sq}}\times(0,1]$ with multiple cylindrical voids removed and inclusions embedded, producing a perforated structure of considerable geometric complexity.
Notably, two removed disks nearly touch at the point $(1.935,1.78)$ (see Figure~\ref{fig:complex_geometric_COMSOL_failure}), 
which is reported to cause mesh generation failures in commercial FEM software such as COMSOL.

\begin{figure}[!htbp]
    \renewcommand\figurename{Figure}
      \centering
      \includegraphics[scale=0.22]{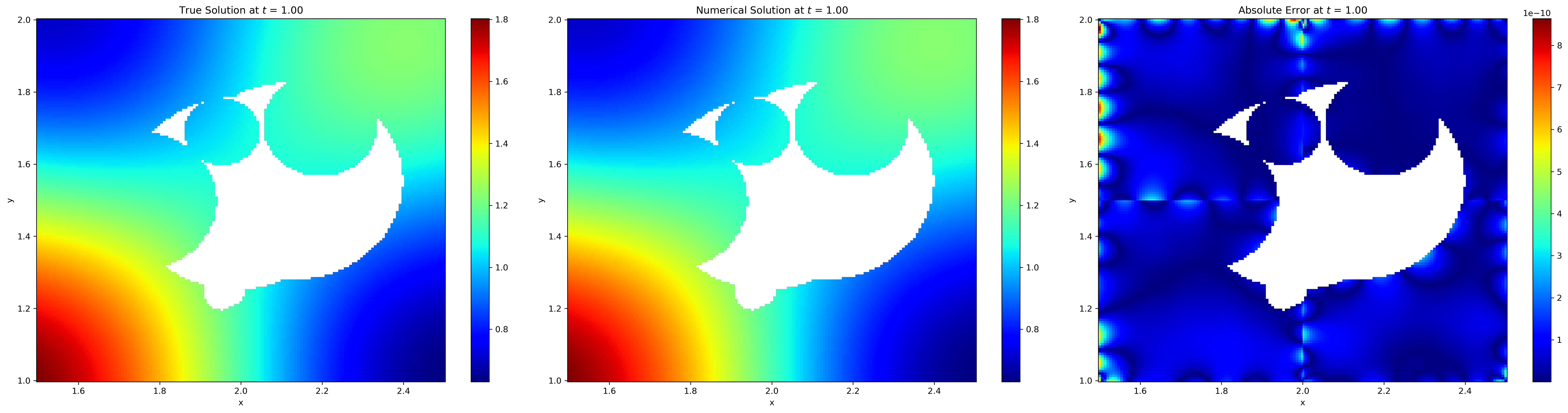}
    \caption{True solution (left), numerical solution (middle), and absolute error (right) at the final time $t=1.00$ for the nonlinear diffusion equation \eqref{eq:complex_diffusion}.}
    \label{fig:complex_diffusion}
\end{figure}

Table~\ref{2D_diffusion_reaction_PDE_result} presents the numerical results of the AMIPN method for \eqref{eq:complex_diffusion} on this geometrically intricate domain. 
As the collocation density $\sqrt[3]{Q}$ and the number of random features $J$ increase, the accuracy improves significantly. 
The relative $L^2$ error of $u$ decreases from $2.01\times10^{-5}$ (at $\sqrt[3]{Q}=20, J=400$) to $6.00\times10^{-10}$ (at $\sqrt[3]{Q}=30, J=2400$), 
while the errors of its first derivatives $(u_x,u_y,u_t)$ decrease from $O(10^{-3}\!-\!10^{-4})$ to $O(10^{-7}\!-\!10^{-8})$. 
Figure~\ref{fig:complex_diffusion} shows that the cross-sectional error at the final time $t=1.0$ reaches the order of $10^{-10}$, 
confirming the capability of the proposed method to deliver high accuracy on geometrically intricate space-time domains where traditional mesh-based methods may fail.

\begin{table}[!htbp]
  \centering
  \caption{The results of the AMIPN method for the nonlinear diffusion problem \eqref{eq:complex_diffusion}}
  \label{2D_diffusion_reaction_PDE_result}
  \begin{tabular}{c c c c c c c c c c}
    \toprule
    $(N_x, N_y, N_z)$ & $\sqrt[3]{Q} $ & $J$ & Solver & IT & NJ & $u$ error & $u_x$ error &$u_y$ error &$u_t$ error \\
    \midrule
    \multirow{7}{*}{$(2, 2, 2)$} 
        & \multirow{2}{*}{$20$} & \multirow{1}{*}{400} & AMIPN & $5$ & $4$& $ 2.01e-05$ & $3.79e-04$ &$3.93e-04$ &$ 2.83e-03$\\
        && \multirow{1}{*}{800} & AMIPN & $5$ & $4$& $ 6.67e-07$ & $1.52e-05$ &$1.61e-05$ &$1.62e-04$\\
    \cmidrule(lr){2-10}
        & \multirow{2}{*}{$25$} 
        &\multirow{1}{*}{1200} & AMIPN  & $5$ & $4$& $ 7.43e-08$ & $2.00e-06$ &$2.08e-06$ &$ 2.14e-05$\\
        && \multirow{1}{*}{1600} & AMIPN & $5$ & $4$& $9.84e-09$ & $2.72e-07$ &$2.89e-07$ &$3.80e-06$\\
      \cmidrule(lr){2-10}
        & \multirow{2}{*}{$30$} 
        & \multirow{1}{*}{2000} & AMIPN & $5$ & $4$& $2.68e-09$ & $7.68e-08$ &$7.73e-08$ &$1.18e-06$\\
        && \multirow{1}{*}{2400} & AMIPN & $5$ & $4$& $6.00e-10$ & $ 1.88e-08$ &$1.95e-08$ &$2.98e-07$\\
    \bottomrule
  \end{tabular}
\end{table}

In summary, the four benchmark problems with moving holes, topological changes, large deformations, and complex space-time geometries demonstrate that the AMIPN method consistently delivers robustness, scalability, and spectral-like accuracy. 
With only a few outer iterations and Jacobian evaluations, errors are reduced to $O(10^{-9}\!-\!10^{-11})$, underscoring its strong potential for solving challenging time-dependent PDEs on domains where conventional mesh-based methods often break down.

\section{Conclusions}\label{sec:7}
In this work, we developed two nonlinear solvers designed for large-scale and ill-conditioned nonlinear least-squares problems arising from RFM discretizations of nonlinear PDEs. The inexact Newton method with right preconditioning  integrates randomized preconditioning with Krylov iterations and derivative-free line searches, thereby improving robustness and stability. Its adaptive multi-step extension   further reduces redundant computations by reusing preconditioned Jacobians with an adaptive early-stopping strategy, resulting in enhanced efficiency without loss of accuracy.

Extensive numerical experiments on three-dimensional steady-state PDEs-including elliptic, diffusion-reaction, Helmholtz, and Gray-Scott equations-as well as two-dimensional timedependent problems on domains with moving holes, topological changes, large deformations, and complex geometries confirm the effectiveness of the proposed solvers. The results consistently show spectral-like convergence, with relative errors reduced to the order of $10^{-9}-10^{-11}$ within only a few outer iterations and a small number of Jacobian evaluations. Compared with conventional FEM and FDM, the proposed methods achieve improvements of four to seven orders of magnitude under the same collocation density, while also outperforming recent machine-learning-based approaches such as PINNs and WANs.
Future work will focus on a rigorous convergence analysis of IPN and AMIPN and on extending the framework to more challenging classes of nonlinear PDEs, further strengthening its applicability and reliability in large-scale scientific and engineering computations.

\appendix
\section{Supplementary numerical examples}\label{Appendix_A}
This appendix presents several additional numerical examples beyond the main test cases reported in the paper.

We first consider a three-dimensional steady-state nonlinear elliptic problem defined on the unit cube $(0,1)^3$:
\begin{equation} \label{3D_nonanalysis_PDE}
\begin{cases}
-\nabla \cdot(\Lambda(u) \nabla u)=1, & (x, y, z) \in(0,1)^3, \\
u=0, & (x, y, z) \in \partial(0,1)^3,
\end{cases}
\end{equation}
where the nonlinear diffusion coefficient is given by $\Lambda(u)=1+u^2$.
Since no closed-form exact solution is available, this problem is used to examine the numerical self-convergence of the proposed method.

As shown in Table~\ref{3D_nonanalysis_nonlinear_elliptic_equation}, 
the proposed AMIPN method exhibits a clear trend of numerical convergence for this 3D nonlinear elliptic problem. 
The solution obtained with the largest parameter setting $(J=2400)$ is taken as the reference solution, 
and the resulting errors for other cases decrease steadily as $J$ increases. 
In particular, when $J=2000$, the error of $u$ reaches the order of $10^{-5}$ ($7.38\times10^{-5}$), 
while the errors of its first derivatives $(u_x,u_y,u_z)$ are on the order of $10^{-4}$ ($8.70\times10^{-4}\sim9.30\times10^{-4}$).

\begin{table}[!htbp]
  \centering
  \caption{The results of the AMIPN method for 3D nonlinear elliptic equation \eqref{3D_nonanalysis_PDE}}
  \label{3D_nonanalysis_nonlinear_elliptic_equation}
  \begin{tabular}{c c c c c c c c c c c}
    \toprule
    $(N_x, N_y, N_z)$ & $\sqrt[3]{Q}$ & $J$ & Solver & IT & NJ &CPU &$u$ error & $u_x$ error & $u_y$ error & $u_z$ error \\
    \midrule
    \multirow{7}{*}{$(2, 2, 2)$} 
        & \multirow{7}{*}{$40$} 
        & $400$ & AMIPN & $2$ & $1$& $114.80$ &$2.37e-03$ & $1.58e-02$ & $1.57e-02$ & $1.61e-02$ \\
        \cmidrule(lr){3-11}
        & & $800$ & AMIPN & $2$ & $1$& $247.13$& $6.60e-04$ & $6.32e-03$ & $6.25e-03$ & $6.31e-03$ \\
        \cmidrule(lr){3-11}
        && $1200$ & AMIPN & $2$ & $1$& $405.85$&$3.09e-04$ & $3.40e-03$ & $3.36e-03$ & $3.39e-03$ \\
        \cmidrule(lr){3-11}
        && $1600$ & AMIPN & $2$ & $1$&$803.16$& $1.60e-04$ & $1.88e-03$ & $1.86e-03$ & $1.87e-03$ \\
        \cmidrule(lr){3-11}
        && $2000$ & AMIPN & $3$ & $1$&$1346.66$& $7.38e-05$ & $8.70e-04$ & $9.30e-04$ & $8.79e-04$ \\
        \cmidrule(lr){3-11}
        & & $2400$ & AMIPN & $3$  & $1$ & $1625.84$& \multicolumn{4}{c}{reference} \\
    \bottomrule
  \end{tabular}
\end{table}

We next consider the two-dimensional Korteweg-de Vries (KdV) initial-boundary value problem defined on the space-time domain $(x, y, t) \in \Omega=[0,1]^2 \times[0,1]$ :
\begin{equation}\label{Kdv_equation}
\begin{cases}
u_t - u\,u_x - u\,u_y + u_{xxx} + u_{yyy} = f(x,y,t), 
 & (x,y,t)\in (0,1)^2\times(0,1], \\[2mm]
u(0,y,t) = g_1(y,t),\quad 
u(1,y,t) = g_2(y,t),\quad 
u_x(0,y,t) = g_3(y,t), 
 & y\in[0,1],\ t\in[0,1], \\[2mm]
u(x,0,t) = p_1(x,t),\quad 
u(x,1,t) = p_2(x,t),\quad 
u_y(x,0,t) = p_3(x,t),
 & x\in[0,1],\ t\in[0,1], \\[2mm]
u(x,y,0) = h(x,y), 
 & (x,y)\in[0,1]^2 .
\end{cases}
\end{equation}
The source term $f$ and the boundary-initial data $\left(g_i, p_i, h\right)$ are constructed from the manufactured exact solution
\begin{equation}
    u(x, y, t)=10-\left(x^3+y^3+t^3\right)
\end{equation}
As shown in Table~\ref{Kdv_equation_results}, the AMIPN method exhibits a clear trend of numerical convergence for the KdV equation. As the collocation density $\sqrt[3]{Q}$ and the number of random features $J$ increase, the relative $L^2$-error of $u$ decreases from $5.48 \times 10^{-6}$ to $1.22 \times 10^{-11}$, while the errors of its derivatives ( $u_x, u_y, u_t$ ) are simultaneously reduced from the order of $10^{-4} \sim 10^{-3}$ to the order of $10^{-9} \sim 10^{-8}$.
Throughout all tested cases, convergence is achieved within only $4-5$ outer iterations, each requiring merely 3 Jacobian evaluations, demonstrating the high efficiency and robustness of the proposed method even for this challenging nonlinear problem.

\begin{table}[!htbp]
  \centering
  \caption{The results of the AMIPN method for the KdV equation \eqref{Kdv_equation}}
  \label{Kdv_equation_results}
  \begin{tabular}{c c c c c c c c c c c}
    \toprule
    $(N_x, N_y, N_t)$ & $\sqrt[3]{Q} $ & $J$ & Solver & IT& NJ & CPU & $u$ error & $u_x$ error & $u_y$ error & $u_t$ error\\
    \midrule
    \multirow{7}{*}{$(2, 2, 2)$} 
        & \multirow{2}{*}{$20$} & \multirow{1}{*}{400} & AMIPN & $4$ & $3$ & $ 50.21$ & $ 5.48e-06$ & $3.68e-04$ & $3.79e-04$ & $1.82e-03$\\
        && \multirow{1}{*}{800} & AMIPN & $4$ & $3$ & $ 139.99$ & $ 6.74e-08$ &$ 6.05e-06$ & $5.83e-06$ & $3.63e-05$\\
    \cmidrule(lr){2-11}
        & \multirow{2}{*}{$25$} 
        &\multirow{1}{*}{1200} & AMIPN & $4$ & $3$ & $  423.17$ & $3.10e-09$ &$3.28e-07$ & $3.26e-07$ & $2.22e-06$\\
        && \multirow{1}{*}{1600} & AMIPN & $4$ & $3$ & $668.82$ & $ 3.35e-10$ &$3.05e-08$ & $3.03e-08$ & $3.06e-07$\\
      \cmidrule(lr){2-11}
        & \multirow{3}{*}{$30$} 
        & \multirow{1}{*}{1600} & AMIPN & $5$ & $3$ & $ 918.26$ & $4.30e-10$ &$3.85e-08$ & $3.78e-08$ & $3.86e-07$\\
        && \multirow{1}{*}{2000} & AMIPN & $5$ & $3$ & $ 1381.28$ & $1.21e-10$ &$5.78e-09$ & $5.66e-09$ & $1.15e-07$\\
        && \multirow{1}{*}{2400} & AMIPN & $5$ & $3$ & $1906.82$ & $1.22e-11$ &$8.68e-10$ & $9.22e-10$ & $1.46e-08$\\
    \bottomrule
  \end{tabular}
\end{table}

We then consider the two-dimensional nonlinear Schr{\"o}dinger equation with Dirichlet boundary conditions:
\begin{equation}\label{SchrodingerEquation}
\left\{
\begin{aligned}
&i\,\frac{\partial h}{\partial t}
+\frac{1}{2}\Big(\frac{\partial^2 h}{\partial x^2}
                 +\frac{\partial^2 h}{\partial y^2}\Big)
+ |h|^2 h
= f(x,y,t), && (x,y,t)\in \ (0, 1)^2 \times(0,1], \\[3pt]
&h(a,y,t)=g_1(y,t),\quad h(b,y,t)=g_2(y,t), && y\in[0,1],\ t\in[0,1],\\
&h(x,c,t)=g_3(x,t),\quad h(x,d,t)=g_4(x,t), && x\in[0,1],\ t\in[0,1],\\
&h(x,y,0)=g_0(x,y), && (x,y)\in[0, 1]^2 .
\end{aligned}
\right.
\end{equation}
The exact solution is given by $h=u+i v$, 
\begin{align*}
    u=(x-1)^3+(y-1)^3+(t-1)^3, \quad v=2 \sin (x(1-x)) \sin (y(1-y)) e^t
\end{align*}
The results in Table~\ref{SchrodingerEquation_results} indicate that the AMIPN method maintains consistently high accuracy across all tested configurations.
Even at moderate settings ( $\sqrt[3]{Q}=20, J=800$ ), the errors are already reduced to $8.79 \times 10^{-7}$ for $u, 6.67 \times 10^{-6}$ for $v$, and $1.22 \times 10^{-6}$ for $h$.
Further increasing $\sqrt[3]{Q}$ and $J$ yields steadily improved accuracy, ultimately reaching $3.04 \times 10^{-10}$, $2.05 \times 10^{-9}$, and $4.01 \times 10^{-10}$ for $u, v$, and $h$, respectively.

\begin{table}[!htbp]
  \centering
  \caption{The results of the AMIPN method for the Schr{\"o}dinger equation \eqref{SchrodingerEquation}}
  \label{SchrodingerEquation_results}
  \begin{tabular}{c c c c c c c c c}
    \toprule
    $(N_x, N_y, N_t)$ & $\sqrt[3]{Q} $ & $J$ & Solver & IT& NJ & $u$ error & $v$ error & $h$ error \\
    \midrule
    \multirow{7}{*}{$(2, 2, 2)$} 
        & \multirow{2}{*}{$20$} & \multirow{1}{*}{400} & AMIPN & $4$ & $3$ & $2.75e-05$ & $ 1.86e-04$ & $3.62e-05$ \\
        && \multirow{1}{*}{800} & AMIPN &$4$ & $3$ & $ 8.79e-07$ & $6.67e-06$ & $1.22e-06$ \\
    \cmidrule(lr){2-9}
        & \multirow{2}{*}{$25$} 
        &\multirow{1}{*}{1200} & AMIPN &$4$ & $3$ & $1.95e-08$ & $1.47e-07$ & $2.70e-08$ \\
        && \multirow{1}{*}{1600} & AMIPN &$4$ & $3$ & $3.18e-09$ & $2.37e-08$ & $4.38e-09$ \\
      \cmidrule(lr){2-9}
        & \multirow{3}{*}{$30$} 
        & \multirow{1}{*}{1600} & AMIPN &$5$ & $3$ & $1.74e-09$ & $1.32e-08$ & $2.42e-09$ \\
        && \multirow{1}{*}{2000} & AMIPN &$4$ & $3$ & $3.35e-10$ & $ 2.34e-09$ & $4.48e-10$ \\
         && \multirow{1}{*}{2400} & AMIPN &$4$ & $3$ & $3.04e-10$ & $ 2.05e-09$ & $4.01e-10$ \\
    \bottomrule
  \end{tabular}
\end{table}

We further solve the following two-dimensional coupled unsteady Burgers equations using the proposed method:
\begin{equation}\label{2D_Burgers}
\begin{aligned}
\partial_t u+u \partial_x u+v \partial_y u&=\frac{1}{\operatorname{Re}}\left(\partial_x^2 u+\partial_y^2 u\right), \\
\partial_t v+u \partial_x v+v \partial_y v&=\frac{1}{\operatorname{Re}}\left(\partial_x^2 v+\partial_y^2 v\right),
\end{aligned}
\qquad x, y \in[0,1],\; t\in[0,1].
\end{equation}
Here, $u$ and $v$ denote the velocity components along the $x$ - and $y$-directions, respectively. The Reynolds number $\operatorname{Re}$ is defined as $\operatorname{Re}=U L / \nu$, where $U$ and $L$ are the characteristic velocity and length, and $\nu$ is the kinematic viscosity of the fluid. The exact solution can be obtained \cite{baeza2006adaptive} as
\begin{align*}
& u(x, y, t)=\frac{3}{4}-\frac{1}{4[1+\exp ((-4 x+4 y-t) \operatorname{Re} / 32)]} \\
& v(x, y, t)=\frac{3}{4}+\frac{1}{4[1+\exp ((-4 x+4 y-t) \operatorname{Re} / 32)]}
\end{align*}
subject to Dirichlet boundary conditions on all boundaries. In this study, the Reynolds number is set to $\operatorname{Re}=100$.

\begin{figure}[!htbp]
    \renewcommand\figurename{Figure}
        \centering
        \includegraphics[scale=0.30]{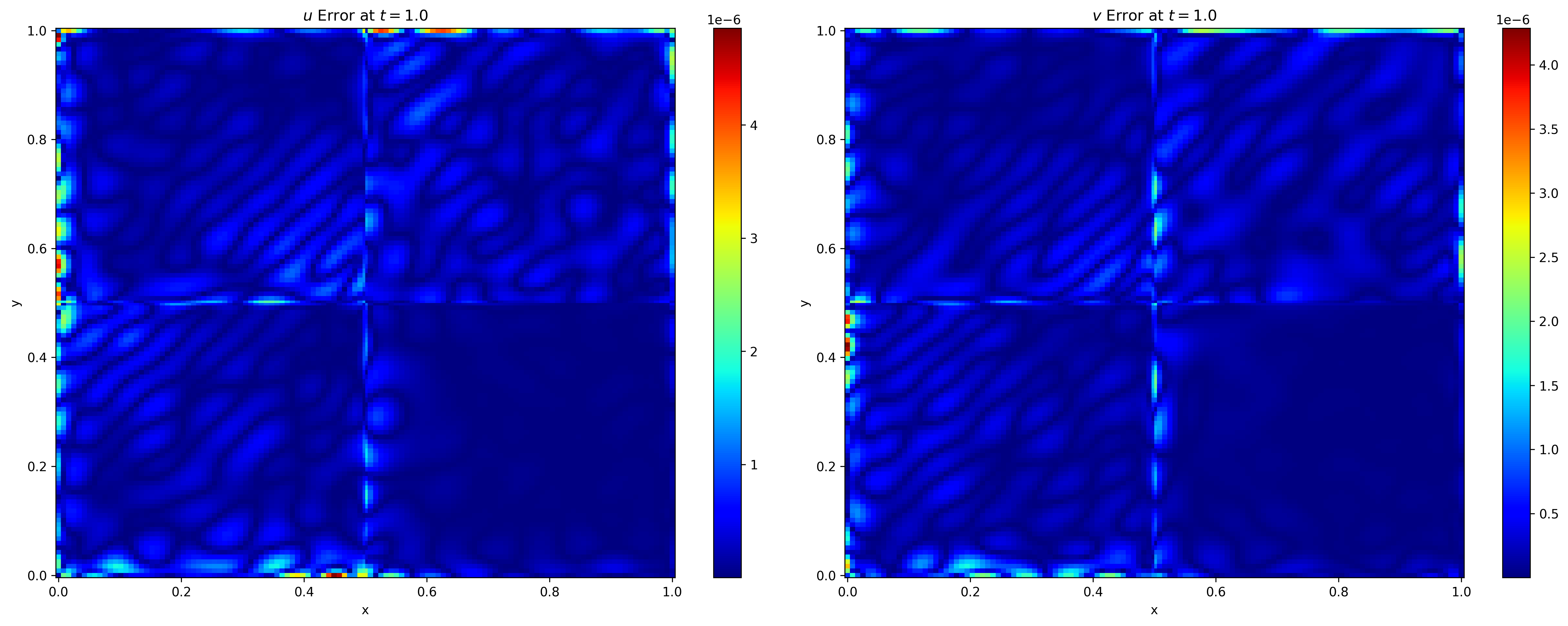}
    \caption{Absolute errors of $u$ (left) and $v$ (right) at the final time $t=1.00$ for the Burgers equations \eqref{2D_Burgers}.}
      \label{fig:2D_Burgers_equation}
\end{figure}

Table~\ref{2D_Burgers_results} presents the numerical results for the Burgers equations. 
As shown in Figure~\ref{fig:2D_Burgers_equation}, the absolute error distributions of $u$ (left) and $v$ (right) at the final time $t=1.00$ 
indicate that the relatively larger errors are localized near the left boundary, 
while the errors remain small in the interior region and along the other boundaries.

\begin{table}[!htbp]
  \centering
  \caption{The results of the AMIPN method for the Burgers equations \eqref{2D_Burgers}}
  \label{2D_Burgers_results}
  \begin{tabular}{c c c c c c c c}
    \toprule
    $(N_x, N_y, N_t)$ & $\sqrt[3]{Q} $ & $J$ & Method & IT & NJ & $u$ error & $v$ error \\
    \midrule
    \multirow{5}{*}{$(2, 2, 2)$} 
        & \multirow{2}{*}{$20$} & \multirow{1}{*}{400} & AMIPN & $4$ &$3$&$4.93e-04$ & $4.15e-04$ \\
        && \multirow{1}{*}{800} & AMIPN & $4$ &$3$&$ 5.11e-05$ & $3.20e-05$ \\
    \cmidrule(lr){2-8}
        & \multirow{2}{*}{$25$} 
        &\multirow{1}{*}{1200} & AMIPN & $4$ &$3$ & $1.15e-05$ & $ 7.56e-06$ \\
        && \multirow{1}{*}{1600} & AMIPN & $4$ &$3$& $3.50e-06$ & $ 2.23e-06$ \\
      \cmidrule(lr){2-8}
        & \multirow{3}{*}{$30$} 
        & \multirow{1}{*}{1600} & AMIPN & $4$ &$3$& $3.87e-06$ & $2.80e-06$ \\
        && \multirow{1}{*}{2000} & AMIPN & $4$ &$3$& $1.17e-06$ & $ 7.20e-07$ \\
        && \multirow{1}{*}{2400} & AMIPN & $4$ &$3$& $6.23e-07$ & $ 3.78e-07$ \\
    \bottomrule
  \end{tabular}
\end{table}

We finally consider the two-dimensional incompressible Navier-Stokes equations to further verify the accuracy of the proposed method. Let the velocity field be $\mathbf{u}=(u(x, y, t), v(x, y, t))$ and the pressure field be $p(x, y, t)$. The governing system is given by
\begin{equation}\label{eq:NS}
\begin{cases}
\frac{\partial u}{\partial t}
+u\frac{\partial u}{\partial x}
+v\frac{\partial u}{\partial y}
=-\frac{\partial p}{\partial x}
+\left(\frac{\partial^2 u}{\partial x^2}
+\frac{\partial^2 u}{\partial y^2}\right),
\\[6pt]
\frac{\partial v}{\partial t}
+u\frac{\partial v}{\partial x}
+v\frac{\partial v}{\partial y}
=-\frac{\partial p}{\partial y}
+\left(\frac{\partial^2 v}{\partial x^2}
+\frac{\partial^2 v}{\partial y^2}\right),
\\[6pt]
\frac{\partial u}{\partial x}
+\frac{\partial v}{\partial y}=0,
\end{cases}
\qquad (x,y,t)\in[0,1]^2\times[0,1].
\end{equation}
Following \cite{ethier1994exact}, the exact solution is prescribed as
\begin{align*}
u(x,y,t)&=-\cos (\pi x)\sin (\pi y)e^{-2\pi^2 t},\\
v(x,y,t)&= \sin (\pi x)\cos (\pi y)e^{-2\pi^2 t},\\
p(x,y,t)&=-\frac{\cos (2\pi x)+\cos (2\pi y)}{4}e^{-4\pi^2 t},
\end{align*}
and Dirichlet boundary conditions are imposed on all boundaries.

Table~\ref{2D_unsteady_N_S} summarizes the numerical results obtained by the AMIPN method, including the number of outer iterations, the number of Jacobian evaluations, and the $L^2$-relative errors of $u$, $v$, and $p$.
As shown in the table, the problem size (i.e., the numbers of equations $m$ and unknowns $n$) grows rapidly with the increase of the collocation density $\sqrt[3]{Q}$ and the number of random basis functions $J$, reaching up to $(745{,}200,48{,}000)$ in the largest case.
Even for such extremely large-scale NLS systems, the AMIPN method achieves rapid convergence within only 4 outer iterations and 3 Jacobian evaluations, while the relative errors decrease steadily—from $10^{-3}\sim 10^{-2}$ at $J=400$ to $10^{-8}\sim10^{-6}$ at $J=2000$—demonstrating both the scalability and the accuracy of the proposed approach for the incompressible Navier-Stokes problems.

\begin{table}[!htbp]
  \centering
  \caption{The results of the AMIPN method for the Navier-Stokes equations \eqref{eq:NS}.}
  \label{2D_unsteady_N_S}
  \begin{tabular}{c c c c c c c c c c}
    \toprule
    $(N_x, N_y, N_t)$ & $\sqrt[3]{Q} $ & $J$ &$(m, n)$& Method & IT& NJ & $u$ error & $v$ error & $p$ error \\
    \midrule
    \multirow{7}{*}{$(2, 2, 2)$} 
        & \multirow{2}{*}{$20$} 
        &\multirow{1}{*}{400} &$(235200, 9600)$& AMIPN & $4$ & $3$  & $ 1.63e-03$ & $1.55e-03$ & $1.63e-02$\\
        \cmidrule(lr){4-10}
        && \multirow{1}{*}{800} & $(235200, 19200)$&AMIPN & $4$ & $3$  & $5.96e-05$ & $ 6.11e-05$ & $8.76e-04$\\
    \cmidrule(lr){2-10}
        & \multirow{2}{*}{$25$} 
        &\multirow{1}{*}{1200} & $(442500, 28800)$&AMIPN & $4$ & $3$  & $ 4.47e-06$ & $ 4.63e-06$ & $1.32e-04$\\
        \cmidrule(lr){4-10}
        && \multirow{1}{*}{1600} & $(442500, 38400)$&AMIPN & $4$ & $3$  & $ 5.43e-07$ & $2.04e-04$ & $3.51e-03$\\
      \cmidrule(lr){2-10}
        & \multirow{2}{*}{$30$} 
        & \multirow{1}{*}{1600} & $(745200, 38400)$&AMIPN & $4$ & $3$  & $5.94e-07$ & $ 6.16e-07$ & $1.44e-05$\\
        \cmidrule(lr){4-10}
        && \multirow{1}{*}{2000} & $(745200, 48000)$&AMIPN & $4$ & $3$  & $9.47e-08$ & $ 9.01e-08$ & $5.65e-06$\\
    \bottomrule
  \end{tabular}
\end{table}

\section*{Acknowledgements}
This work is partially supported by the National Key R\&D Program of China (No. 2022YFA1005200, No. 2022YFA1005202, and No.2022YFA1005203), the NSFC Major Research Plan - Interpretable and General Purpose Next-generation Artificial Intelligence (Nos. 92270001 and 92370205), NSFC grant 12425113, and Key Laboratory of the Ministry of Education for Mathematical Foundations and Applications of Digital Technology, University of Science and Technology of China.

\section*{CRediT authorship contribution statement}
Longze Tan: Writing-review \& editing, Writing-original draft, Visualization, Validation, Software, Methodology, Investigation, Data curation, Mathematical analysis, Formal analysis, Conceptualization. Xiaohe Yue: Validation, Software, Visualization, Data curation. Jingrun Chen: Writing-review \& editing, Writing-original draft, Visualization, Validation, Supervision, Project administration, Methodology, Investigation, Funding acquisition, Formal analysis, Conceptualization. Jiamin Jiang: Writing-review \& editing, Writing-original draft, Data curation, Methodology, Formal analysis. 

 \section*{Declaration of competing interest}
 The authors confirm that they have no financial or personal conflicts of interest that could have influenced the findings presented in this paper.

 \section*{Data availability}
 No data was used for the research described in the article.

\bibliographystyle{elsarticle-num}
\bibliography{bibliography}

\end{document}